%% file: cnlohner.tex
\documentclass[10pt]{article}
\usepackage{amssymb}
\usepackage{amsmath}
\usepackage{graphicx}

\def\oarrow{{\mathrel{\longrightarrow\mkern-25mu\circ}\;\;}}
\newcommand{\Id}{\mathrm{Id}}
\def\qed{{\hfill{\vrule height5pt width3pt depth0pt}\medskip}}

\newfont{\bbc}{msbm10}
\def\Bbb#1{\hbox{{\bbc #1}}}
\newcommand{\inte }{{\rm int}\,}

\newcommand{\dom }{\mathrm{dom} \,}
\newcommand{\ip}[2]{\left\langle#1|#2\right\rangle}
\newcommand{\orbit}[1]{\mathcal O\left(#1\right)}
\newcommand{\Fix}{\mathrm{Fix}}

\newtheorem{theorem}{Theorem}
\newtheorem{lemma}[theorem]{Lemma}
\newtheorem{remark}[theorem]{Remark}
\newtheorem{definition}[theorem]{Definition}
\begin{document}
\begin{center}
    {\bf \LARGE $\mathcal C^r$-Lohner algorithm}
    \vskip\baselineskip
    {\large
    Daniel Wilczak\footnote{
            Research supported by an annual national scholarship for young scientists from
            the Foundation for Polish Science
        },
    Piotr Zgliczy\'nski\footnote{
            Research supported in part by Polish State Ministry of
            Science and Information Technology  grant N201 024 31/2163 
         }
    } \\
 Jagiellonian University, Institute of Computer Science, \\
 Nawojki 11, 30--072  Krak\'ow, Poland \\ e-mail: wilczak@ii.uj.edu.pl, umzglicz@cyf-kr.edu.pl
\vskip\baselineskip \today
\end{center}
\begin{abstract}
We present a Lohner type algorithm  for the computation of
rigorous bounds for solutions of ordinary differential equations
and its derivatives with respect to initial conditions up to
arbitrary order. As an application we prove the existence of
multiple invariant tori around some elliptic periodic
orbits for the pendulum equation with periodic forcing and for
Michelson system.
\end{abstract}

\input intro.tex

\input notation.tex

\input equation.tex
\input cnalg.tex

\input enclosure.tex
\input eval.tex
\input poincare.tex
\input pendulum.tex
\input msystem.tex
\input implementation.tex
\input appendix.tex

\end{document}

%% file: intro.tex
\section{Introduction}
This paper is a sequel to \cite{Z}. We present here a Lohner-type
algorithm for computation of rigorous enclosures of  partial
derivatives with respect to initial conditions  up to an arbitrary
order $r$ of the flow induced by an autonomous ODE, hence the name
$\mathcal C^r$-Lohner algorithm. Let $r$ be a positive integer,
then  by $\mathcal C^r$-algorithm  we will mean the routine which gives
rigorous estimates for partial derivatives with respect to initial
conditions up to an order $r$ and $\mathcal C^r$-computations we mean an
application of an $\mathcal C^r$-algorithm.

Our main  motivation for the development of $\mathcal C^r$-algorithm was a
desire to provide a tool, which will considerably extend the
possibilities of computer assisted proofs in the dynamics of ODEs.
Till now most of such proofs have used  topological conditions
(see for example \cite{HZHT,MM,GZ,Z1}) and additionally conditions
on the first derivatives with respect to initial conditions (see
for example \cite{RNS,T,Wi1,WZ,KZ}), hence it required $\mathcal C^0$- and
$\mathcal C^1$-computations, respectively. The spectrum of problems treated
 includes the questions of the existence of periodic orbits and
their local uniqueness, the existence of symbolic dynamics, the
existence of hyperbolic invariants sets, the existence of homo-
and heteroclinic orbits. To treat other phenomena, like
bifurcations of periodic orbits, the route to chaos, invariant
tori through KAM theory one needs the knowledge of partial
derivatives with respect to initial conditions of higher order.

In principle, one can think that a good rigorous ODE solver should
be enough. Namely, to compute the partial derivatives of the flow
induced by
\begin{equation}
   x'=f(x), \qquad x \in \mathbb{R}^n \label{eq:odeeq}
\end{equation}
 it is enough to rigorously integrate a system of variational
 equations obtained by  a formal differentiation of
 (\ref{eq:odeeq}) with respect to the initial conditions. For
 example for $r=2$ we have the following system
\begin{eqnarray}
  x'&=&f(x),  \label{eq:var0}\\
  \frac{d}{dt} V_{ij}(t) &=& \sum_{s=1}^n \frac{\partial f_i}{\partial
  x_s}(x)  V_{sj}(t)  \label{eq:var1} \\
  \frac{d}{dt} H_{ijk}(t) &=& \sum_{s,r=1}^n \frac{\partial^2
  f_i}{\partial x_s \partial x_r
  }(x)V_{rk}(t)V_{sj}(t)+
    \sum_{s=1}^n \frac{\partial f_i}{\partial x_s}(x)
  H_{sjk}(x), \label{eq:var2}
\end{eqnarray}
with the initial conditions
\begin{eqnarray*}
   x(0)=x_0, \quad V(0)=Id, \qquad H_{ijk}(0)=0, \quad
   i,j,k=1,\dots,n. \label{eq:varic}
\end{eqnarray*}
It is well known that if by  $\varphi(t,x_0)$ we denote the
(local) flow induced by (\ref{eq:odeeq}), then
\begin{eqnarray*}
  \frac{\partial \varphi_i}{\partial x_j}(t,x_0) &=& V_{i,j}(t), \\
 \frac{\partial^2 \varphi_i}{\partial x_j \partial
  x_k}(t,x_0) &=&  H_{ijk}(t).
\end{eqnarray*}
Analogous statements are true for higher order partial derivatives
with respect to initial conditions.

It turns out that a straightforward application of a rigorous ODE
solver to the system of variational Equations
(\ref{eq:var0}--\ref{eq:var2}) is  very inefficient. Namely, it
totally ignores the structure of the system and leads to a very
poor performance and unnecessary long computation times (see
Section~\ref{subsec:whycr}).

Our algorithm is a modification of the Lohner algorithm \cite{Lo},
which takes into account the structure of variational Equations
(\ref{eq:var0}--\ref{eq:var2}). Basically it consists of the
Taylor method, a heuristic routine for a priori bounds for
solution of (\ref{eq:var0}--\ref{eq:var2}) during a time step and
a Lohner-type control of the wrapping effect, which is done
separately for  $x$ and partial derivatives with respect initial
conditions (the variables $V$ and $H$ in
(\ref{eq:var1},\ref{eq:var2})). The Taylor method is realized
using the automatic differentiation \cite{Ra} and the algorithms
for computation of compositions of multivariate Taylor series.

The proposed algorithm  has been successfully applied in
\cite{HNW} to the Michelson system \cite{Mi}, where a computer
assisted proof of the existence of a cocoon bifurcation was
presented. Some parts of this proof  required $\mathcal C^2$-computations.

In the present paper in Section~\ref{sec:pendulum} we show an
application of our algorithm to pendulum equation with periodic
forcing and the Michelson system. We used it to compute rigorous
bounds for the coefficients of some normal forms up to order five,
which enabled us to prove the existence of invariant tori
around some elliptic periodic orbits in these systems
using KAM theorem for twist maps on the plane. These
proofs required $C^3$ and $C^5$ computations.

%% file: notation.tex
\section{Basic definitions}
To effectively deal with the formulas involving  partial
derivatives  we will use extensively a notation of multiindices,
multipointers and submultipointers throughout the paper.

As an motivation  let us consider the formula for the partial
derivatives of the composition of maps. Assume $g:\mathbb R^n\to
\mathbb R^n$ and $f:\mathbb R^n\to\mathbb R$ are of class
$\mathcal C^3$. We have
\begin{eqnarray*}
\frac{\partial^3(f\circ g)}{\partial x_i\partial x_j\partial x_c}
    = \sum_{k,r,s=1}^n
        \frac{\partial^3 f}{\partial x_k\partial x_r\partial x_s}
        \frac{\partial g_k}{\partial x_i}
        \frac{\partial g_r}{\partial x_j}
        \frac{\partial g_s}{\partial x_c} +
    \sum_{k=1}^n
        \frac{\partial f}{\partial x_k}\frac{\partial^3 g_k}{\partial x_i\partial x_j\partial x_c}\\
    +
        \sum_{k,r=1}^n \frac{\partial^2 f}{\partial x_k\partial x_r}\left(
        \frac{\partial^2 g_k}{\partial x_i\partial x_c}\frac{\partial g_r}{\partial x_j}+
        \frac{\partial g_k}{\partial x_i}\frac{\partial^2 g_r}{\partial x_j\partial x_c}+
        \frac{\partial^2 g_k}{\partial x_i\partial x_j}\frac{\partial g_r}{\partial x_c}
        \right)
\end{eqnarray*}
To the operator $\frac{\partial^3}{\partial x_{i_1}\partial
x_{i_2}\partial x_{i_3}}$ we can in a unique way assign a
\emph{multipointer}, which is a nondecreasing sequence of integers
$(j_1,j_2,j_3)$, such that $\{i_1,i_2,i_3\}=\{j_1,j_2,j_3\}$. A
\emph{submultipointer} is  a multipointer, which is a part of a
longer multipointer, for example $(i,j,c)_{(1,3)}=(i,c)$. One
observes, that submultipointers appear at several places in the
above formula.

A \emph{multiindex} is an element of $\alpha\in\mathbb N^n$. It is
another way to represent various partial derivatives.  The
coefficient $\alpha_i$ tells us how many times to differentiate a
function with respect to the $i$-th variable. Obviously, we have
one-to-one correspondence between multipointers and multiindices.

\subsection{Multiindices}
By $\mathbb N$ we will denote the set of nonnegative integers,
i.e. $\mathbb N=\{0,1,2,\ldots\}$.
\begin{definition}
An element $\tau\in\mathbb N^n$ will be called a \emph{multiindex}.
\end{definition}
For a sequence $\alpha=(\alpha_1,\ldots,\alpha_n)\in \mathbb N^n$
and a vector $x=(x_1,\ldots,x_n)\in\mathbb R^n$ we set
\begin{enumerate}
\item $|\alpha|=\alpha_1+\cdots+\alpha_n$
\item $\alpha! = \alpha_1!\cdot\alpha_2!\cdots\alpha_n!$
\item $x^\alpha = (x_1^{\alpha_1},\ldots,x_n^{\alpha_n})$
\end{enumerate}
By $e_i^n\in \mathbb N^n$ we will denote
\begin{equation*}
e_i^n=(0,0,\ldots,0,\overbrace{1}^{i},0,\ldots,0,0).
\end{equation*}
We will drop  the index $n$  (the dimension) in the symbol $e^n_i$
when it  is obvious from the context.

Put $\mathbb N^n_p:=\{a\in\mathbb N^n: |a|=p\}$.

For $\delta=(\delta_1,\ldots,\delta_k)\in\mathbb
N^{n_1}\times\cdots\times\mathbb N^{n_k}$ we set
\begin{enumerate}
\item $|\delta| = \sum_{i=1}^k |\delta_i|$
\item $\delta! = \prod_{i=1}^k \delta_i!$
\end{enumerate}

Let $f=(f_1,\ldots,f_m):\mathbb R^n\to\mathbb R^m$ be sufficiently smooth.
 For $\alpha\in \mathbb N^n$ we set
\begin{enumerate}
\item $\displaystyle D^\alpha f_i = \frac{\partial^{|\alpha|}
f_i}{\partial x_1^{\alpha_1}\cdots\partial x_n^{\alpha_n}}$
\item $D^\alpha f=(D^\alpha f_1,D^\alpha f_2,\ldots,D^\alpha f_m)$
\end{enumerate}
For a function $f:\mathbb R\times\mathbb R^n\to\mathbb R^n$ by
$D^\alpha f_i(t,x)$ we will denote $D^\alpha f_i(t,\cdot)(x)$ and
similarly
\begin{equation*}
D^\alpha f(t,x) = (D^\alpha f_1(t,x),\ldots, D^\alpha f_n(t,x)).
\end{equation*}
This convention means that $D^\alpha$ always acts on $x$-variables.

\subsection{Multipointers}
For a fixed $n>0$ and $p>0$ we define
\begin{eqnarray*}
    \mathcal N_p^n&:=&\left\{(a_1,a_2,\ldots,a_p) \in \mathbb{N}^p : 1\leq a_1\leq\cdots\leq a_p\leq n\right\}\\
    \mathcal N=\mathcal N^n&:=&\bigcup_{p=1}^\infty\mathcal N_p^n
\end{eqnarray*}
\begin{definition}
An element of $\mathcal N^n$ will be called a \emph{multipointer}.
\end{definition}
\begin{remark}
A function
\begin{equation}\label{eq:defLambda}
\Lambda:\mathcal N^n_p\ni (a_1,\ldots,a_p)\to \sum_{i=1}^p e^n_{a_i}
\in\mathbb N^n_p
\end{equation}
is a bijection.
\end{remark}
Let $f=(f_1,\ldots,f_m):\mathbb R^n\to\mathbb R^m$ be a sufficiently smooth.
 For $a\in \mathcal
N_p^n$ we set
\begin{enumerate}
\item $\displaystyle D_a f_i := \frac{\partial^p f_i}{\partial x_{a_1}\ldots\partial x_{a_p}}$
\item $\displaystyle D_a f := (D_a f_1,\ldots,D_a f_m)$
\end{enumerate}
For a function $f:\mathbb R\times\mathbb R^n\to\mathbb R^n$ by
$D_a f_i(t,x)$ we will denote $D_af_i(t,\cdot)(x)$. In the light
of the above notations $D_\alpha f = D^{\Lambda(\alpha)}f$.

For $a=(a_1,a_2,\ldots,a_n)\in\mathbb N^n_p$ and
$b=(b_1,b_2,\ldots,b_n)\in\mathbb N^n_q$ we define
\begin{equation*}
a+b=(a_1+b_1,\ldots,a_n+b_n)\in\mathbb N^n_{p+q}.
\end{equation*}
For $\alpha\in\mathcal N^n_p$ and $\beta\in\mathcal N^n_q$ we
define
\begin{equation*}
\alpha+\beta=
\Lambda^{-1}\left(\Lambda(\alpha)+\Lambda(\beta)\right)\in\mathcal
N^n_{p+q}.
\end{equation*}

By $\leq$ we will denote a linear order (lexicographical
order) in $\mathcal N$ defined in the following way. For
$a\in\mathcal N^n_p$ and $b\in\mathcal N^n_q$
\begin{equation}\label{eq:linearOrder}
\left(a\leq b\right) \Longleftrightarrow
    \begin{cases}
    \text{either } \exists i, i\leq p, i \leq q, a_i<b_i \text{ and }  a_j=b_j \text{ for } j<i\\
    \text{or } p\leq q \text{ and } a_i=b_i \text{ for }
    i=1,\ldots,p.
    \end{cases}
\end{equation}

\begin{definition}
For $k\leq p$ we set
\begin{equation}\label{eq:defNpk}
\mathcal N^p(k):=\{(\delta_1,\ldots,\delta_k)\in(\mathcal N^p)^k :
\delta_1\leq\cdots\leq\delta_k,
\delta_1+\cdots+\delta_k=(1,2,\ldots,p)\}
\end{equation}
\end{definition}

We will use $\mathcal N^p(k)$ extensively in the next section. Its
will be used to label terms in $D^\alpha f_i(\varphi(t,x))$.
Observe that for $p>0$
\begin{eqnarray*}
  \mathcal{N}^p(1)=\{ (1,2,\dots,p) \} \\
  \mathcal{N}^p(p)=\{ ((1),(2),\dots,(p)) \}
\end{eqnarray*}
One can construct all elements of $\mathcal{N}^p(k)$ using the
following recursive procedure. From the definition of $\mathcal
N^p(k)$  it follows that if
$(\delta_1,\ldots,\delta_{m-1})\in\mathcal N^{p-1}(m-1)$ then
$(\delta_1,\ldots,\delta_{m-1},(p))\in\mathcal N^p(m)$ (notice
that order is preserved). Similarly, if
$(\delta_1,\ldots,\delta_m)\in\mathcal N^{p-1}(m)$ then
$$(\delta_1,\ldots,\delta_{s-1},\delta_s+(p),\delta_{s+1},\ldots,\delta_m)\in\mathcal
N^p(m)$$ and again order of elements is preserved. Hence, for
$p>2$ and $1<k<p$ we have $\mathcal N^p(k) = A\cup B$ where
\begin{equation}\label{eq:indexDecomposition}
\begin{split}
A &= \left\{
        (\delta_1,\ldots,\delta_{k-1},(p)) : (\delta_1,\ldots,\delta_{k-1})\in\mathcal N^{p-1}(k-1)
    \right\}\\
B&=
    \bigcup_{s=1}^k\left\{
    (\delta_1,\ldots,\delta_{s-1},\delta_s+(p),\delta_{s+1},\ldots,\delta_k) :(\delta_1,\ldots,\delta_k)\in\mathcal N^{p-1}(k)
    \right\}
\end{split}
\end{equation}
and the sets $A$ and $B$ are disjoint.

Another way to generate all elements of $\mathcal{N}^p(k)$ can be
described as follows
\begin{itemize}
\item decompose the set $\{1,2,\dots,p\}$ into $k$ nonempty and
disjoints sets $\Delta_i$, $i=1,\dots,k$
\item we sort each $\Delta_i$ and permute $\Delta_i$'s to obtain $\min(\Delta_1) < \min(\Delta_2) < \dots <
  \min(\Delta_k)$
\item we define $\delta_i$ to be an ordered set consisting of all
elements of $\Delta_i$ for $i=1,\dots,k$
\end{itemize}

\begin{definition}
For an arbitrary $a\in\mathcal N_p^n$ and $\delta\in\mathcal
N^p_k$ such that $k\leq p$ we define a \emph{submultipointer}
$a_\delta \in\mathcal N^n_k$ by $(a_\delta)_i = a_{\delta_i}$ for
$i=1,\ldots,k$, which can be expressed using $\Lambda$ as follows
\begin{equation*}
a_\delta := \Lambda^{-1}\left(\sum_{i=1}^k
e^n_{a_{\delta_i}}\right)\in\mathcal N^n_k
\end{equation*}
\end{definition}

%% file: equation.tex
\section{Equations for variations}
\label{sec:eqvar}

Consider an ODE $x'=f(x)$ where $f$ is $\mathcal C^{K+1}$. Let
$\varphi:\mathbb R\times\mathbb R^n\oarrow \mathbb R^n$ be a local
dynamical system induced by $x'=f(x)$. It is well known, that
$\varphi \in \mathcal C^{K}$ and one can derive  the equations for
partial derivatives of $\varphi$
 by differentiating equation
 $\frac{\partial \varphi}{\partial t}(t,x)=f(\varphi(t,x))$
 with respect to the initial condition $x$.
As a result we obtain a system of so-called equations for
variations, whose size depends on the order $r$ of partial
derivatives we intend to compute. An example of such system for
$r=2$ is given by (\ref{eq:var0}--\ref{eq:var2}) with initial
conditions given by (\ref{eq:varic}).

The goal of this section is to write the equations for variations
in a compact form using multipointers and  multiindices, which
allows us to take into account the symmetries of partial
derivatives,

\begin{lemma}
Assume $f\in\mathcal C^{r+1}$ and let $\varphi:\mathbb R\times \mathbb
R^n\oarrow\mathbb R^n$ be a local dynamical system induced by
$x'=f(x)$. Then for $a\in\mathcal N^n_p$ such that $p\leq r$ holds
\begin{equation}\label{eq:eqDer}
\frac{d}{dt}D_a \varphi_i
    =   \sum_{k=1}^{p}\
        \sum_{i_1,\ldots,i_k=1}^n
        \left(D^{e_{i_1}+\cdots+e_{i_k}}
        f_i\right)\circ\varphi
        \sum_{(\delta_1,\ldots,\delta_k)\in \mathcal N^p(k)}\
        \prod_{j=1}^k D_{a_{\delta_j}}\varphi_{i_j}
\end{equation}
for $i=1,\ldots,n$.
\end{lemma}

\textbf{Proof:} In the proof the functions
$D^{e_{i_1}+\cdots+e_{i_k}} f_i$ are always evaluated at
$\varphi(t,x)$, and various partial derivatives of $\varphi$ are
always evaluated at $(t,x)$, therefore the arguments will be
always dropped to simplify formulae. We prove the lemma by
induction on $p=|a|$. If $p=1$ then $a=(c)$ for some
$c\in\{1,\ldots,n\}$ and \eqref{eq:eqDer} becomes
\begin{equation*}
   \frac{d}{dt} D_{(c)}\varphi_i= \frac{d}{dt}\frac{\partial\varphi_i}{\partial x_c}
    = \sum_{s=1}^n \frac{\partial f_i}{\partial x_s}\frac{\partial \varphi_s}{\partial x_c}
    =\sum_{s=1}^n D^{e_{s}}f_i\cdot D_{(c)}\varphi_{s}.
\end{equation*}

Assume \eqref{eq:eqDer} holds true for $p-1$, $p>1$. Let us fix
$a\in\mathcal N^n_p$. We have $a=b+(c)$, where
$b=(a_1,\ldots,a_{p-1})\in\mathcal N^n_{p-1}$ and $c=a_p$. Since
\eqref{eq:eqDer} is satisfied for $p-1$, therefore  we have
\begin{equation*}
\begin{aligned}
\frac{d}{dt}D_a \varphi_i &=
D_{(c)}\left(\frac{d}{dt}D_b\varphi_i\right)\\
    &=   D_{(c)}\left(
        \sum_{k=1}^{p-1}\
        \sum_{\substack{i_1,\ldots,i_k=1\\ \beta:=e_{i_1}+\cdots+e_{i_k}}}^n
        D^{\beta} f_i
        \sum_{(\delta_1,\ldots,\delta_k)\in \mathcal N^{p-1}(k)}\
        \prod_{j=1}^k D_{b_{\delta_j}}\varphi_{i_j}
    \right)\\
    &=
        \sum_{k=1}^{p-1}\
        \sum_{\substack{i_1,\ldots,i_{k+1}=1\\
        \beta:=e_{i_1}+\cdots+e_{i_{k+1}}}}^n
        D^{\beta} f_i\cdot D_{(c)}\varphi_{i_{k+1}}
        \sum_{(\delta_1,\ldots,\delta_k)\in \mathcal N^{p-1}(k)}\
        \prod_{j=1}^k D_{b_{\delta_j}}\varphi_{i_j}
            \\
            &+
        \sum_{k=1}^{p-1}\
        \sum_{\substack{i_1,\ldots,i_k=1\\ \beta:=e_{i_1}+\cdots+e_{i_k}}}^n
        D^{\beta} f_i
        \sum_{(\delta_1,\ldots,\delta_k)\in \mathcal N^{p-1}(k)}\
            \sum_{s=1}^k D_{b_{\delta_s}+(c)}\varphi_{i_s}\prod_{\substack{j=1,\\ j\neq s}}^k D_{b_{\delta_j}}\varphi_{i_j}
    \end{aligned}
\end{equation*}

For $k=1,\dots,p$ we set
\begin{equation}
  T_k:= \sum_{i_1,\ldots,i_k=1}^n
        D^{e_{i_1}+\cdots+e_{i_k}} f_i
        \sum_{(\delta_1,\ldots,\delta_k)\in \mathcal N^p(k)}\
        \prod_{j=1}^k D_{a_{\delta_j}}\varphi_{i_j}
\end{equation}
Now our goal is to prove that:
\begin{equation}
\frac{d}{dt}D_a \varphi_i  = \sum_{k=1}^{p}  T_k
\end{equation}
Our strategy of proof is as follows. We will define $S_1,\dots,S_p$,
such that
\begin{equation}
  \frac{d}{dt}D_a \varphi_i = \sum_{k=1}^p S_k, \qquad
  S_i=T_i, \quad i=1,\dots,p.
\end{equation}
We set
\begin{equation*}
\begin{aligned}
  S_1&=\sum_{k=1}\sum_{\substack{i_1,\ldots,i_k=1\\ \beta:=e_{i_1}+\cdots+e_{i_k}}}^n
        D^{\beta} f_i
        \sum_{(\delta_1,\ldots,\delta_k)\in \mathcal N^{p-1}(k)}\
   \sum_{s=1}^k D_{b_{\delta_s}+(c)}\varphi_{i_s}\prod_{\substack{j=1,\\ j\neq s}}^k
   D_{b_{\delta_j}}\varphi_{i_j}\\
  S_p&= \sum_{k=p-1}\
        \sum_{\substack{i_1,\ldots,i_{k+1}=1\\ \beta:=e_{i_1}+\cdots+e_{i_{k+1}}}}^n
        D^{\beta} f_i\cdot D_{(c)}\varphi_{i_{k+1}}
        \sum_{(\delta_1,\ldots,\delta_k)\in \mathcal N^{p-1}(k)}\
        \prod_{j=1}^k D_{b_{\delta_j}}\varphi_{i_j}.
\end{aligned}
\end{equation*}
For $m=2,3,\ldots,p-1$ we set
\begin{equation*}
\begin{aligned}
S_m &=
        \sum_{k=m-1}\
        \sum_{\substack{i_1,\ldots,i_{k+1}=1\\ \beta:=e_{i_1}+\cdots+e_{i_{k+1}}}}^n
        D^{\beta} f_i\cdot D_{(c)}\varphi_{i_{k+1}}
        \sum_{(\delta_1,\ldots,\delta_k)\in \mathcal N^{p-1}(k)}\
        \prod_{j=1}^k D_{b_{\delta_j}}\varphi_{i_j}
            \\
            &+
        \sum_{k=m}\
        \sum_{\substack{i_1,\ldots,i_k=1\\ \beta:=e_{i_1}+\cdots+e_{i_k}}}^n
        D^{\beta} f_i
        \sum_{(\delta_1,\ldots,\delta_k)\in \mathcal N^{p-1}(k)}\
            \sum_{s=1}^k D_{b_{\delta_s}+(c)}\varphi_{i_s}\prod_{\substack{j=1,\\ j\neq s}}^k D_{b_{\delta_j}}\varphi_{i_j}
\end{aligned}
\end{equation*}

It remains to show that $S_i=T_i$ for $i=1,\dots,p$. Consider first
$i=1$. Recall that $\mathcal{N}^{p-1}(1)=\{(1,2,\dots,p-1)\}$, hence
\begin{eqnarray*}
    S_1=\sum_{s=1}^n D^{e_s}f_i \cdot D_{b+(c)}\varphi_{s}
    = \sum_{s=1}^n D^{e_s}f_i \cdot D_a\varphi_{s}.
\end{eqnarray*}
Therefore
\begin{equation}
  S_1=T_1.
\end{equation}

Consider now $i=p$. For an arbitrary $s>0$ $\mathcal N^s(s)$
contains only one element $((1),(2),\ldots,(s))$. Therefore we
obtain
\begin{equation*}
\begin{aligned}
S_{p}&= \sum_{i_1,\ldots,i_p=1}^n
        D^{e_{i_1}+\cdots+e_{i_p}} f_i\cdot D_{(c)}\varphi_{i_p}
        \sum_{(\delta_1,\ldots,\delta_{p-1})\in \mathcal N^{p-1}(p-1)}\
        \prod_{j=1}^{p-1} D_{b_{\delta_j}}\varphi_{i_j} \\
&=
        \sum_{i_1,\ldots,i_p=1}^n
        D^{e_{i_1}+\cdots+e_{i_p}} f_i\cdot D_{(c)}\varphi_{i_p}
        \prod_{j=1}^{p-1} D_{b_j}\varphi_{i_j}.
\end{aligned}
\end{equation*}
Since $a=b+(c)$, where $c=(a_p)$, hence
\begin{equation*}
\begin{aligned}
S_{p}&=
        \sum_{i_1,\ldots,i_p=1}^n
        D^{e_{i_1}+\cdots+e_{i_p}} f_i\prod_{j=1}^{p}
        D_{a_j}\varphi_{i_j}\\
        &=
        \sum_{i_1,\ldots,i_p=1}^n
        D^{e_{i_1}+\cdots+e_{i_p}} f_i
        \sum_{(\delta_1,\ldots,\delta_p)\in \mathcal N^p(p)}\
        \prod_{j=1}^{p} D_{a_{\delta_j}}\varphi_{i_j}=T_p
\end{aligned}
\end{equation*}

Consider now $m=2,3,\ldots,p-1$. We have
\begin{equation*}
\begin{aligned}
 S_m&= \sum_{i_1,\ldots,i_m=1}^n
        D^{e_{i_1}+\cdots+e_{i_m}} f_i\cdot D_{(c)}\varphi_{i_m}
        \sum_{(\delta_1,\ldots,\delta_{m-1})\in \mathcal N^{p-1}(m-1)}\
        \prod_{j=1}^{m-1} D_{b_{\delta_j}}\varphi_{i_j}
            \\
            &+
        \sum_{i_1,\ldots,i_m=1}^n
        D^{e_{i_1}+\cdots+e_{i_m}} f_i
        \sum_{(\delta_1,\ldots,\delta_m)\in \mathcal N^{p-1}(m)}\
            \sum_{s=1}^m D_{b_{\delta_s}+(c)}\varphi_{i_s}\prod_{\substack{j=1,\\ j\neq s}}^m D_{b_{\delta_j}}\varphi_{i_j}
\end{aligned}
\end{equation*}
Using decomposition $\mathcal N^p(m)=A\cup B$ as in
\eqref{eq:indexDecomposition} we obtain
\begin{equation*}
\begin{aligned}
S_m&=
        \sum_{i_1,\ldots,i_m=1}^n
        D^{e_{i_1}+\cdots+e_{i_m}} f_i
        \sum_{(\delta_1,\ldots,\delta_{m-1},\delta_m=(p))\in A}\
        \prod_{j=1}^m D_{a_{\delta_j}}\varphi_{i_j}
            \\
       &+
        \sum_{i_1,\ldots,i_m=1}^n
        D^{e_{i_1}+\cdots+e_{i_m}} f_i
        \sum_{(\delta_1,\ldots,\delta_m)\in B}\
            \prod_{j=1}^m D_{a_{\delta_j}}\varphi_{i_j}\\
    &=
        \sum_{i_1,\ldots,i_m=1}^n
        D^{e_{i_1}+\cdots+e_{i_m}} f_i
        \sum_{(\delta_1,\ldots,\delta_m)\in \mathcal N^p(m)}\
        \prod_{j=1}^m D_{a_{\delta_j}}\varphi_{i_j}=T_m
\end{aligned}
\end{equation*}
We have shown that $T_i=S_i$ for $i=1,\dots,p$. This finishes the
proof. \qed

%% file: cnalg.tex
\section{$\mathcal C^r$-Lohner algorithm}
\label{sec:cralg}

\subsection{Why one needs an $\mathcal C^r$-algorithm?}
\label{subsec:whycr}
 There are several effective algorithms for the computation of rigorous
bounds for solutions of ordinary differential equations, including
Lohner method \cite{Lo}, Hermite--O\-bresch\-koff algorithm
\cite{NJ} or Taylor models \cite{BM}.  For $\mathcal
C^r$-com\-pu\-ta\-tions the number of
equations to solve is equal to $n\begin{pmatrix}n+r\\
n \end{pmatrix}$  hence, even for $r=1$ direct application of such
an algorithms to equations for variations \eqref{eq:mainEquation}
leads to integration in high dimensional space and is usually
inefficient. Let us recall after \cite[Sec. 6]{Z}  the basic
reason for this. In order to have a good control over the
expansion rate of the set of initial conditions during a time step
 these algorithms, while being $\mathcal C^0$, are $\mathcal C^1$ 'internally'(or
higher for Taylor models), because they solve non-rigorously
equations for ($\frac{\partial \varphi}{\partial x}$) - the
variational matrix of the flow. This effectively squares the
dimension of phase space of the equation and impacts heavily the
computation time. But as it was observed in \cite{Z} the equations
for partial derivatives of the flow can be seen as non-autonomous
and nonhomogenous linear equations, therefore we do not need
additional equations for variations for them. As a result the
dimension of the effective phase space for our
$\mathcal C^r$-algorithm is given by $n\begin{pmatrix}n+r\\
n \end{pmatrix}$ and not a square of this number.

Another important aspect of the proposed algorithm is the fact that
the Lohner-type control of the wrapping effect is done separately
for $x$-variables and  variables  $D_a \varphi$. This feature is not
present in the blind application of $\mathcal C^0$ algorithm to the
system of variational equations and it turns out that this often
practically switches off the control of the wrapping effect on
$x$-variables, as various choices used in this control become
dominated by the $D_a \varphi$-variables.

 In \cite{Z} a $\mathcal C^1$-algorithm has been proposed. Here we present
 an algorithm for computation of higher order partial derivatives.

\subsection{An outline of the algorithm}

Let us fix $r\leq K$ and consider the following system of
differential equations
\begin{equation}\label{eq:mainEquation}
\left\{
\begin{aligned}
    \frac{d}{dt}\varphi
    &=  f\circ \varphi\\
        \frac{d}{dt}D_a \varphi
    &=  \sum_{k=1}^{d}\
        \sum_{i_1,\ldots,i_k=1}^n
        \left(D^{e_{i_1}+\cdots+e_{i_k}} f\right)\circ\varphi
           \sum_{(\delta_1,\ldots,\delta_k)\in \mathcal N^{d}(k)}\
        \prod_{j=1}^k D_{a_{\delta_j}}\varphi_{i_j}
\end{aligned}
\right.
\end{equation}
for all $a\in\mathcal N^n_d$, $d=1,\ldots,r$.

Our goal is to present an algorithm for computing a rigorous bound
for the solution of \eqref{eq:mainEquation} with a set of initial
conditions
\begin{equation}\label{eq:initCondition}
\begin{cases}
    \varphi(0,x_0)&\in [x_0]\subset\mathbb R^n\\
    D\varphi(0,x_0)&=\mathrm{Id}\\
    D_a\varphi(0,x_0)&=0,\qquad \text{for } a\in\mathcal N^n_2\cup\ldots\cup\mathcal
    N^n_r.
\end{cases}
\end{equation}

In the sequel we will use the following notations:
\begin{itemize}
\item if a solution of system \eqref{eq:mainEquation} is defined for $t>0$ and some
$x_0\in\mathbb R^n$, then for $a \in \mathcal{N}$  by $V_a(t,x_0)$
we denote $D_a\varphi(t,x_0)$
\item for $[x_0]\subset \mathbb R^n$ by $[V_a(t,[x_0])]$ we will denote a set for
which we have $V_a(t,[x_0])\subset [V_a(t,[x_0])]$. This set  is
obtained using an rigorous numerical routine described below.
\end{itemize}

The $\mathcal C^r$-Lohner algorithm is a modification of $\mathcal
C^1$-Lohner algorithm \cite{Z}. One step of $\mathcal C^r$-Lohner is
a shift along the trajectory of the system \eqref{eq:mainEquation}
with the following input and output data
\newline {\bf Input data:}
\begin{itemize}
\item $t_k$ - a current time,
\item $h_k$ - a time step,
\item $[x_k] \subset {\Bbb R}^n $,
    such that $\varphi(t_k,[x_0]) \subset [x_k]$,
\item $[V_{k,a}]=[V_{k,a}(t_k,[x_0])]\ \subset {\Bbb R}^n$, such that  $D_a\varphi(t_k,[x_0]) \subset
[V_{k,a}]$ for $a\in\mathcal N^n_1\cup\ldots\cup\mathcal N^n_r$.
\end{itemize}

\noindent {\bf Output data:}
\begin{itemize}
\item $t_{k+1}=t_{k} + h_k$ - a new current time,
\item $[x_{k+1}] \subset {\Bbb R}^n $,  such that $\varphi(t_{k+1},[x_0]) \subset
[x_{k+1}]$,
\item $[V_{k+1,a}]=[V_{k+1,a}(t_{k+1},[x_0])]\ \subset {\Bbb R}^n$,
such that  $D_a\varphi(t_{k+1},[x_0]) \subset
[V_{k+1,a}]$ for $a\in\mathcal N^n_1\cup\ldots\cup\mathcal N^n_r$.
\end{itemize}
We will often skip the arguments of $V_{k,a}$ when they are
obvious from the context.

The values of $[x_{k+1}]$ and $[V_{k+1,a}]$, $a\in\mathcal N^n_1$
are computed using one step $\mathcal C^1$-Lohner algorithm. After
it is done, we perform the following operations  to compute
$[V_{k+1,a}]$ for $a\in\mathcal N^n_2\cup\ldots\cup\mathcal N^n_r$
\begin{description}
\item[1.] Find a rough enclosure for $D_a\varphi([0,h_k],[x_k])$.
\item[2.] Compute $[V_{k+1,a}]$, this will also involve some
rearrangement computations to reduce the wrapping effect for $V$
\cite{Mo,Lo}.
\end{description}

%% file: enclosure.tex
\section{Computation of a rough enclosure for $D_a\varphi$}
\label{sec:encllognorm} For a fixed multipointer $a\in\mathcal
N^n_d$ Equation \eqref{eq:mainEquation} can be  written as follows
\begin{equation}\label{eq:AB}
  \frac{d}{dt}D_a\varphi(t,x)=B_a(t,x) + A(t,x) D_a\varphi(t,x)
\end{equation}
where
\begin{equation}\label{eq:defAB}
\begin{aligned}
 B_a&= \sum_{k=2}^{d}\
        \sum_{i_1,\ldots,i_k=1}^n
        \left(D^{e_{i_1}+\cdots+e_{i_k}} f\right)\circ\varphi
           \sum_{(\delta_1,\ldots,\delta_k)\in \mathcal N^{d}(k)}\
        \prod_{j=1}^k D_{a_{\delta_j}}\varphi_{i_j}\\
 A&= Df\circ\varphi
\end{aligned}
\end{equation}

The procedure for computing the rough enclosure is based on the
notion of a logarithmic norm, which we give below.
\begin{definition}\cite{HNW}
For a square matrix $A$ the logarithmic norm $\mu(A)$ is defined as
a limit
\begin{equation*}
\mu(A)=\limsup_{h\to 0^+}\frac{\left\|\Id+Ah\right\|-1}{h}
\end{equation*}
where $\|\cdot\|$ is a given matrix norm.
\end{definition}

The  formulas for the logarithmic norm of a real matrix in the
most frequently used norms are (see \cite{HNW})
\begin{enumerate}
\item for $\|x\|_1=\sum_i |x_i|$,  $\mu(A)=\max_{j}(a_{jj}+\sum_{i\neq j}|a_{ij}|)$
\item for m $\|x\|_2=\sqrt{\sum_i |x_i|^2}$,  $\mu(A)$ is equal to the
largest eigenvalue of $(A+A^T)/2$
\item for  $\|x\|_\infty=\max_i |x_i|$,  $\mu(A)=\max_{i}(a_{ii}+\sum_{j\neq i}|a_{ij}|)$
\end{enumerate}

 In order to find bounds for $D_a\varphi$ we use the following
theorem \cite[Thm. I.10.6]{HNW}
\begin{theorem}
\label{thm:LnHNW} Let $x(t)$  be a solution of a  differential
equation
\begin{equation}
  x'(t)=f(t,x(t)), \quad x \in {\Bbb R}^n
\end{equation}
Let $\nu(t)$ be a piecewise differentiable function with values in
${\Bbb R}^n$. Assume that
\begin{eqnarray*}
  \mu\left(\frac{\partial f}{\partial x}(t,\eta)\right) \leq l(t)
  \quad \mbox{for  $\eta \in [x(t),\nu(t)]$ } \\
  |\nu'(t) - f(t,\nu(t))| \leq \delta(t),
\end{eqnarray*}
where by $\mu(A)$,  we denote a logarithmic norm of a square matrix
$A \in {\Bbb R}^{n \times n}$.

 Then for $t \geq t_0$ we have
\begin{equation}
  |x(t) - \nu(t)| \leq e^{L(t)}\left(|x(t_0) - \nu(t_0)| + \int_{t_0}^t e^{-L(s)}\delta(s) ds
  \right),
\end{equation}
with $L(t)=\int_{t_0}^t l(\tau) d\tau$.
\end{theorem}

We apply the above theorem to Equation (\ref{eq:AB}) to obtain
\begin{lemma}
Let us fix $x\in\Bbb R^n$. Assume that $|B_a(t,x)| \leq
\delta(t) $ and $\mu(A(t,x)) \leq l(t)$, then for $t>t_0$
\begin{equation}
  |D_a\varphi(t,x)| \leq |D_a\varphi(t_0,x)| e^{L(t)}+
  e^{L(t)}\int_{t_0}^te^{-L(\tau)} \delta(\tau)d\tau
\end{equation}
with $L(t)=\int_{t_0}^t l(\tau) d\tau$.
\end{lemma}
\textbf{Proof:} Consider Equation \eqref{eq:AB} and a homogenous
problem for (\ref{eq:AB})
\begin{equation}
\frac{d}{dt}w=  f(t,w):=A(t,x) \cdot w, \qquad w \in {\Bbb R}^n.
\label{eq:encllin}
\end{equation}
Using Theorem \ref{thm:LnHNW} we can estimate the difference
between any solution of (\ref{eq:encllin}), $w$, and a solution of
\eqref{eq:AB}, denoted by $D_a\varphi$.
\begin{equation}
  |D_a\varphi(t) - w(t)| \leq |D_a\varphi(t_0) - w(t_0)| e^{L(t)}+
  e^{L(t)}\int_{t_0}^te^{-L(\tau)} \delta(\tau)d\tau.
\end{equation}
After a substitution $w(t)=0$, which is a solution of the
homogenous equation, we obtain our assertion. \qed

Usually, we do not have any control over the time dependence of
$\delta$ and $l$, hence we will use the following
\begin{lemma}
Assume that $|B_a(t,x)| \leq \delta$  and $\mu(A(t,x)) \leq l$ for
$t \in [0,h]$ then for $t \in [0,h]$ we have
\begin{equation}
  |D_a\varphi(t,x)| \leq |D_a\varphi(0,x)|\max(1, e^{hl})+
  \delta \frac{e^{lt}-1}{l}, \quad \mbox{if $l \neq 0$},
\end{equation}
or
\begin{equation}
  |D_a\varphi(t,x)| \leq |D_a\varphi(0,x)|+
  \delta t, \quad \mbox{when $l = 0$}.
\end{equation}
\end{lemma}

\subsection{The procedure for the computation of the rough enclosure for $V$.}
The procedure for the computing of the rough enclosure is
iterative, which means that given a rough enclosure for
$\varphi([0,h_k],[x_k])$ and rough enclosures
$D_a\varphi([0,h_k],[x_k])$ for all $a\in\mathcal
N^n_1\cup\ldots\cup \mathcal N^n_p$ we are able to compute the
rough enclosure for $D_a\varphi([0,h_k],[x_k])$ for $a\in\mathcal
N^n_{p+1}$.

The procedures for computation of the rough enclosures of
$\varphi([0,h_k],[x_k])$ and $D_a\varphi([0,h_k],[x_k])$ for
$a\in\mathcal N^n_1$ has been given in \cite{Z}. Below we present
an algorithm for computing $[E_a]$ for $a\in\mathcal
N^n_2\cup\ldots\cup\mathcal N^n_r$.

\noindent {\bf Input parameters:}
\begin{itemize}
\item $h_k$ - a time step,
\item $[x_k] \subset {\Bbb R}^n $ - the current value of
$x=\varphi(t_k,[x_0])$,
\item $[E_0]\subset {\Bbb R}^n$ - a compact and convex  such that
$\varphi([0,h_k],[x_k])\subset [E_0]$
\item $[E_{a}]\subset {\Bbb R}^n$, $a\in\mathcal
N^n_1\cup\ldots\cup\mathcal N^n_p$ such that
$D_a\varphi([0,h_k],[x_k])\subset [E_a]$ for $a\in\mathcal
N^n_1\cup\ldots\cup\mathcal N^n_p$.
\end{itemize}

\noindent {\bf Output:}
\begin{itemize}
\item $[E_a] \subset {\Bbb R}^n$, $a\in\mathcal N^n_{p+1}$ such that
\begin{equation*}
  D_a\varphi([0,h_k],[x_k])\subset [E_a]
\end{equation*}
\end{itemize}

Before we present an algorithm let us observe that for a fixed
$a\in\mathcal N^n_{p+1}$, $B_a$ defined in \eqref{eq:defAB} could
be seen as a multivariate function of $t$, $x$ and
$V_b=D_b\varphi$ for $b\in\mathcal N^n_1\cup\ldots\cup\mathcal
N^n_p$. More precisely, put $m_p:=\sharp\left\{\mathcal
N^n_1\cup\ldots\cup\mathcal N^n_p\right\}$, where $\sharp$ stands
for number of elements of a set. Recall that,  we have defined by
\eqref{eq:linearOrder}  a linear order  in $\mathcal N^n$. Hence,
there is a unique sequence of multipointers $b_1,\ldots,b_{m_p}$,
such that $b_i\in\mathcal N^n_1\cup\ldots\cup\mathcal N^n_p$ for
$i=1,\ldots,m_p$, $b_1\leq b_2\leq\cdots\leq b_{m_p}$ and $b_i\neq
b_j$ for $i\neq j$.

Let us define
\begin{eqnarray*}
    \tilde B_a:\mathbb R\times\left(\mathbb R^n\right)^{m_p+1}\rightarrow\mathbb R^n,\\
    F_a:\mathbb R\times\left(\mathbb R^n\right)^{m_p+1}\rightarrow\mathbb R^n
\end{eqnarray*}
by
\begin{equation}
\label{eq:defBTilde}
\begin{split}
\tilde B_a(t,x,v_{b_1},\ldots,v_{b_{m_p}}) =
        \sum_{k=2}^{p+1}\
        \sum_{i_1,\ldots,i_k=1}^n
        D^{e_{i_1}+\cdots+e_{i_k}} f(\varphi(t,x))\\
           \sum_{(\delta_1,\ldots,\delta_k)\in \mathcal N^{p+1}(k)}\
        \prod_{j=1}^k \left(v_{a_{\delta_j}}\right)_{i_j}
\end{split}
\end{equation}
and
\begin{equation}\label{eq:defF}
F_a(t,x,v_{b_1},\ldots,v_{b_m}) = \tilde
B_a(t,x,v_{b_1},\ldots,v_{b_m}) + Df(\varphi(t,x))V_a(t,x)
\end{equation}

\noindent {\bf Algorithm:} \newline To compute $[E_a]$ for
$a\in\mathcal N^n_{p+1}$ we proceed as follows
\begin{description}
\item[1.] Find $l \geq \left(\max_{x \in [E_0]} \mu\left(Df(x)\right)\right)$.
\item[2.] Compute $\delta_a\geq \max \|\tilde B_a\|$, i.e.
\begin{equation*}
 \delta_a\geq \max_{
            \substack{(x,v_{b_1},\ldots,v_{b_{m_p}})\in[E_0]\times[E_{b_1}]\times\cdots\times[E_{b_{m_p}}]}
        }\left\|
            \tilde B_a(0,x,v_{b_1},\ldots,v_{b_{m_p}})
        \right\|
\end{equation*}
For example, if $a=(j,c)\in\mathcal N^n_2$, then $\delta_a$ should
be such that
\begin{equation*}
    \delta_a\geq\max_{x\in[E_0],v_1\in[E_{(1)}],\ldots,v_n\in[E_{(n)}]}
    \left\|
        \sum_{r,s=1}^n\frac{\partial^2 f}{\partial x_r\partial
        x_s}(x)\left(v_j\right)_s\left(v_c\right)_r
    \right\|
\end{equation*}
\item[3.] Define $[E_a]_i=[-1,1] \delta_a \frac{e^{lt}-1}{l}$, for $i=1,\ldots,n$, where $[E_a]_i$ denotes $i$-th coordinate of $[E_a]$.
\end{description}

One can refine the obtained enclosure by
\begin{eqnarray*}
  [E_a]&:=& \left(
    [0,h_k] F_a\left(0,[E_0],[E_{b_1}],\ldots,[E_{b_{m_p}}]\right)
  \right)\cap [E_a]
\end{eqnarray*}
Indeed, for $i=1,\ldots,n$, $t\in[0,h_k]$ and $x_0\in[E_0]$ we
have
\begin{eqnarray*}
D_a\varphi_i(t,x_0) &=& D_a\varphi_i(t,x_0)- D_a\varphi(0,x_0)\\
    &=& t
    \left(F_a\right)_i(\theta_i,x_0,D_{b_1}\varphi(\theta_i,x_0),\ldots,D_{b_{m_p}}\varphi(\theta_i,x_0))\\
    &=& t
    \left(F_a\right)_i(0,\varphi(\theta_i,x_0),D_{b_1}\varphi(\theta_i,x_0),\ldots,D_{b_{m_p}}\varphi(\theta_i,x_0))
\end{eqnarray*}
for some  $\theta_i\in[0,t]\subset[0,h_k]$. In the above we have
used the fact that
\begin{equation*}
F_a(t,x,v_1,\ldots,v_{m_p})=F_a(0,\varphi(t,x),v_1,\ldots,v_{m_p}).
\end{equation*}
Since $\varphi(\theta_i,x_0)\in[E_0]$ and
$D_{b_j}\varphi(\theta_i,x_0)\in[E_{b_j}]$ for $j=1,\ldots,m_p$ we
get
\begin{equation*}
D_a\varphi_i(t,x_0)\in    [0,h_k]
\left(F_a\right)_i\left(0,[E_0],[E_{b_1}],\ldots,[E_{b_{m_p}}]\right)
\end{equation*}

%% file: eval.tex
\section{Computation of $[V_{k+1}]$}
\label{sec:compH}

\subsection{Composition formulas}

For any $p$-times continuously differentiable functions
$f,g:\mathbb{R}^n \to \mathbb{R}^n$ and $a\in\mathcal N^n_{p}$ we
have
\begin{equation}
\label{eq:compf} D_a(f\circ g)= \sum_{k=1}^{p}\
        \sum_{i_1,\ldots,i_k=1}^n
        \left(D^{e_{i_1}+\cdots+e_{i_k}}
        f_i\right)\circ g
        \sum_{(\delta_1,\ldots,\delta_k)\in \mathcal N^p(k)}\
        \prod_{j=1}^k D_{a_{\delta_j}}g_{i_j}
\end{equation}
We can apply the above formula to $f=\varphi(h_k,\cdot)$ and
$g=\varphi(t_k,\cdot)$ to obtain
\begin{eqnarray*}
  V_{a}(t_k+h_k,x_0) = \sum_{k=1}^{p}\
        \sum_{i_1,\ldots,i_k=1}^n
        V_{\Lambda^{-1}(e_{i_1}+\ldots+e_{i_k})}(h_k,x_k)\\
        \sum_{(\delta_1,\ldots,\delta_k)\in \mathcal N^p(k)}\
        \prod_{j=1}^k \left(V_{a_{\delta_j}}\right)_{i_j}(t_k,x_0)
\end{eqnarray*}
for all $x_0\in[x_0]$. Using notations $[V_{k+1,a}]:=[V_a(t_k+h_k,[x_0])]$ and $[V_{k,a}]=[V_a(t_k,[x_0])]$
we can rewrite the above equation as
\begin{equation}
  [V_{k+1,a}] = \sum_{k=1}^{p}\
        \sum_{i_1,\ldots,i_k=1}^n
        V_{\Lambda^{-1}(e_{i_1}+\ldots+e_{i_k})}(h_k,[x_k])\\
        \sum_{(\delta_1,\ldots,\delta_k)\in \mathcal N^p(k)}\
        \prod_{j=1}^k \left[V_{k,a_{\delta_j}}\right]_{i_j}
        \label{eq:compH}
\end{equation}
where $\Lambda$ is defined by \eqref{eq:defLambda}.

\subsection{The procedure for computation of $[V_{k+1}]$}

We introduce  new parameters $o_d$ - the order of the Taylor
method used in computations of $V_a$ for $a\in\mathcal N^n_d$. It
makes sense to take $o_1\geq o_2\geq\cdots\geq o_r$.

\noindent
 {\bf Input parameters:}
 \begin{itemize}
   \item $h_k$ - a time step,
   \item $[x_k] \subset {\Bbb R}^n $ - the current value of
   $x=\varphi(t_k,[x_0])$,
    \item $[V_{k,a}] \subset {\Bbb R}^n$ - a current value of $V_{k,a}(t_k,[x_0])$, for $a\in\mathcal N^n_1\cup\ldots\cup\mathcal N^n_r$
   \item $[E_0] \subset {\Bbb R}^n$ compact and convex, such that
           $\varphi([0,h_k],[x_k]) \subset [E_0]$ - a rough enclosure for
           $[x_k]$,
   \item $[E_a] \subset {\Bbb R}^n$, compact and convex, such that
      $D_a\varphi([0,h_k],[x_k]) \subset [E_a]$, for $a\in\mathcal N^n_1\cup\ldots\cup\mathcal N^n_r$.
 \end{itemize}

\noindent {\bf Output:}
 $[V_{k+1,a}] \subset {\Bbb R}^n$, such that
 \begin{equation}
   V_a(t_k+h_k,x_0)\in [V_{k+1,a}]
 \end{equation}
for $x_0\in[x_0]$ and $a\in\mathcal N^n_1\cup\ldots\cup\mathcal N^n_r$.

\noindent
 {\bf Algorithm:}
We compute $[V_{k+1}]$ as follows
\begin{description}
\item[1.] Computation of $V_a(h_k,[x_k])$ using Taylor method for
Equation \eqref{eq:mainEquation}, i.e. for $a\in\mathcal N^n_p$ we
compute
\begin{eqnarray}
  [F_a]&=&\sum_{i=1}^{o_p} \frac{h_k^i}{i!}
  \frac{d^{i-1}}{dt^{i-1}}F_a(0,[x_k],V_{b_1},\ldots,V_{b_{m_{p-1}}}) \label{eq:Fa} \\
   &+&
  \frac{h^{o_p + 1}}{(o_p +1)!}
  \frac{d^{o_p}}{dt^{o_p}}F_a(0,[E_0],[E_{b_1}],\ldots,[E_{b_{m_{p-1}}}]). \nonumber
\end{eqnarray}
where $V_{b_i}=0$ for $b_i\in\mathcal N^n_2\cup\ldots\cup\mathcal
N^n_{p-1}$ and $V_{(j)}=e^n_j$ for $j=1,\ldots,n$.
  Observe that
  \begin{equation}
     V_a(h_k,[x_k]) \subset [F_a]
  \end{equation}
Indeed, using Taylor series expansion we obtain that for
$x_k\in[x_k]$ and $j=1,\ldots,n$ holds
\begin{eqnarray*}
(V_a)_j(h_k,x_k) = \sum_{i=1}^{o_p} \frac{h_k^i}{i!}
  \frac{d^{i-1}}{dt^{i-1}}(F_a)_j(0,x_k,V_{b_1}(0,x_k),\ldots,V_{b_{m_{p-1}}}(0,x_k))\\
  +\frac{h^{o_p + 1}}{(o_p +1)!}
  \frac{d^{o_p}}{dt^{o_p}}(F_a)_j(\theta_i,x_k,V_{b_1}(\theta_i,x_k),\ldots,V_{b_{m_{p-1}}}(\theta_i,x_k))
\end{eqnarray*}
for some $\theta_i\in[0,h_k]$. Observe, that
\begin{eqnarray*}
\frac{d^{o_p}}{dt^{o_p}}(F_a)_j(\theta_i,x_k,V_{b_1}(\theta_i,x_k),\ldots,V_{b_{m_{p-1}}}(\theta_i,x_k))\\=
\frac{d^{o_p}}{dt^{o_p}}(F_a)_j(0,\varphi(\theta_i,x_k),V_{b_1}(\theta_i,x_k),\ldots,0,V_{b_{m_{p-1}}}(\theta_i,x_k))
\end{eqnarray*}

Using $\varphi(\theta_i,x_k)\in[E_0]$ and
$V_{b_s}(\theta_i,x_k)\in[E_{b_s}]$ for $s=1,\ldots,m_{p-1}$ we
obtain our assertion.
\item[2.] The composition. Put
\begin{equation*}
[J_k]:=([F_{(1)}],\ldots,[F_{(n)}])^T
\end{equation*}
Using \eqref{eq:compH} for $a\in\mathcal N^n_p$ we have
\begin{equation}
  [V_{k+1,a}] =  [\alpha_a] + [J_k]\cdot  [V_{k,a}],   \label{eq:propH}
\end{equation}
where
\begin{equation}
[\alpha_a]=\sum_{k=2}^{p}\
        \sum_{i_1,\ldots,i_k=1}^n
        [F_{\Lambda^{-1}(e_{i_1}+\ldots+e_{i_k})}]
        \sum_{(\delta_1,\ldots,\delta_k)\in \mathcal N^p(k)}\
        \prod_{j=1}^k \left[V_{k,a_{\delta_j}}\right]_{i_j}
\end{equation}
\end{description}

In our implementation of the algorithm we use the symbolic
differentiation to obtain formulae for $D_af$. Next, using the
automatic differentiation we compute
$\frac{d^i}{dt^i}F_a(t,x,V_{b_1}(t,x),\ldots,V_{b_{m_{p-1}}}(t,x))|_{t=0}$
which appear in \eqref{eq:Fa}.

\subsection{Rearrangement for $V_a$ - the evaluation of Equation
(\ref{eq:propH})}

It is well know that a direct evaluation of Equation
(\ref{eq:propH}) leads to wrapping effect \cite{Mo,Lo}.  To avoid it
following the work of Lohner \cite{Lo} we will use the same scheme
as it was proposed in \cite{Z}.

Namely, observe  that  Equation (\ref{eq:propH}) has exactly the
same structure as the propagation equations for $\mathcal
C^1$-method (see \cite[Section 3]{Z}). Moreover, all vectors
$V_{k,a}$, for $a\in\mathcal N^n_1\cup\ldots\mathcal N^n_r$
'propagate' by the same $[J_k]$ as did the variational part in
\cite{Z}, hence it makes sense the same approach.

To be more precise, each set $[V_{k,a}]$, for $a\in\mathcal
N^n_1\cup\ldots\cup\mathcal N^n_r$ is represented in the following
form
\begin{equation*}
[V_{k,a}] = v_{k,a} + [B_{k}][r_{k,a}] + C_k[q_{k,a}]
\end{equation*}
 where $[B_k]$ is interval matrix, $C_k$ is point matrix, $v_{k,a}$
is a point vector and $r_{k,a},q_{k,a}$ are interval vectors.
Observe that  $[B_k]$ and $C_k$ are independent of $a$.

 In
the sequel we will drop index $a$. Equation \eqref{eq:propH} leads
to
\begin{equation}
{}[V_{k+1}]=[\alpha]+[J_k](v_{k} + [B_{k}][r_{k}] + C_k[q_{k}])
  \label{eq:evVk1}
\end{equation}
Let  $m([z])$ denotes a center of an interval object, i.e. $[z]$ is
interval vector or interval matrix and $\Delta([z])=[z]-m([z])$.

Let $[Q]$ be an interval matrix which contains an orthogonal
matrix. Usually, $[Q]$ is computed by the orthonormalisation of
the columns of $m([J_k])[B_k]$.

Let
\begin{eqnarray*}
   \left[Z\right] &=& m([J_k]) C_k\\
   C_{k+1} &=& m([Z])\\
   \left[B_{k+1}\right] &=& [Q]
\end{eqnarray*}

Then we  rearrange formula (\ref{eq:evVk1}) as follows
\begin{equation}\label{eq:evalProcedure}
\begin{array}{rcl}
    \left[s\right] &=& [\alpha]+[J_k]v_k+\Delta([J_k])[V_k]\\
    v_{k+1} &=& m([s])\\
    \left[q_{k+1}\right] &=& [q_k]\\
    \left[r_{k+1}\right] &=& [Q^T](\Delta([s])+\Delta([Z])[q_k])+([Q^T]m([J_k])[B_k])[r_k]
\end{array}
\end{equation}

Summarizing, we can use the following data structure to represent
$\varphi(t_k,[x_0])$ and $D_a\varphi(t_k,[x_0])$, for $a\in\mathcal
N^n_1\cup\ldots\cup\mathcal N^n_r$
\begin{quote}
    \textbf{type} CnSet = \textbf{record}
    \begin{quote}
    $v_0, r_0, q_0$: IntervalVector;\\
    $C_0, B_0, C, B$ : IntervalMatrix;\\
    $\{v_a, r_a, q_a: \text{ IntervalVector}\}_{a\in\mathcal
    N^n_{1}\cup\ldots\cup\mathcal N^n_r}$
    \end{quote}
    \textbf{end;}
\end{quote}
The set $\varphi(t_k,[x_0])$ is represented as $v_0+B_0r_0 +
C_0q_0$, the partial derivatives $D_a\varphi(t_k,[x_0])$ are
represented as $v_a+Br_a + Cq_a$. The matrices $B,C$ are common for
all partial derivatives.

Notice, that if we start the $\mathcal C^r$ computation with an
initial condition \eqref{eq:initCondition} then there is no
Lipschitz part at the beginning for the partial derivatives.
Hence, the initial values for $C$ and $B$ are set to the identity
matrix and the initial values for $q_a,r_a$ are set to zero.

If  the interval vectors $r_a$ become 'thick' (i.e. theirs
diameters are larger than some threshold value) we can set a new
Lipschitz part in our representation (it must be done
simultaneously for all $D_a\varphi$) and reset $r_a$ in the
following way
\begin{eqnarray*}
    q_a &=& r_a + (B^T C)q_a,\quad\text{for }a\in\mathcal N^n_1\cup\ldots\cup\mathcal N^n_r\\
    r_a &=& 0,\quad\text{for }a\in\mathcal N^n_1\cup\ldots\cup\mathcal N^n_r\\
    C &=& B\\
    B &=& \Id
\end{eqnarray*}
A similar change of the Lipshitz part may be done when vectors
$r_a$ become thick in comparison to $q_a$.

%% file: poincare.tex
\section{Derivatives of Poincar\'e map}

Consider a differential equation
\begin{equation}
  x'=f(x), \quad x \in {\Bbb R}^n, \quad f \in \mathcal C^{K+1} \label{eq:ODE}
\end{equation}
Let $\varphi:{\Bbb R} \times {\Bbb R}^n \to {\Bbb R}^n $ be a
(local) dynamical system induced by (\ref{eq:ODE}). Let $\alpha:
{\Bbb R}^n \to {\Bbb R}$ be $\mathcal C^1$-map. Put $\Pi=\{x \ |\
\alpha(x)=C \}$.
\begin{definition}
We will say that $\Pi$ is a local section for the vector field $f$
at $y_0\in\Pi$ if
\begin{equation}
  \ip{\nabla \alpha(y_0)}{f(y_0)} \neq 0. \label{eq:ftransa}
\end{equation}
\end{definition}

Assume $x_0\in \mathbb{R}^n$ and $t_0 \in \mathbb{R}$ are such that
$\Pi$ is a local section at $\varphi(t_0,x_0)$. Consider an implicit
equation
\begin{equation}
  \alpha(\varphi(t_P(x),x))=C.     \label{eq:rowTp}
\end{equation}
It follows easily from (\ref{eq:ftransa}) and from the implicit
function theorem that there exists a uniquely defined $t_P:
\mathbb{R}^n \oarrow \mathbb{R}$  in a neighborhood of $x_0$, such
that $t_P(x_0)=t_0$. The function $t_P$ is as smooth as the flow
$\varphi$. We will refer to $t_P$ as to the \emph{Poincare return
time to section $\Pi$}.

We define a Poincar\'e map  $P:\mathbb{R}^n \supset \dom(t_P) \to
\mathbb{R}^n$ by
\begin{equation}
  P(x)=\varphi(t_P(x),x).  \label{eq:defPoinc}
\end{equation}
Usually the Poincar\'e  map is defined as a map $P:\Pi_1 \oarrow
\Pi_2$, where $\Pi_1,\Pi_2$ are local sections in $\mathbb{R}^n$.
The approach taken here, i.e. treating the Poincar\'e map as map
$P:\mathbb{R}^n \oarrow \mathbb{R}^n$ allows us to not to worry
about the coordinates on local section.

In this section we are interested in the partial derivatives of
$P$ defined by (\ref{eq:defPoinc}).

From (\ref{eq:defPoinc}) we can compute $\frac{\partial
P_i}{\partial x_j}$ and we obtain
\begin{equation}
 \frac{\partial P_i}{\partial x_j}(x) = f_i(P(x)) \frac{\partial t_P}{\partial x_j}(x)  +
  \frac{\partial \varphi_i}{\partial x_j}(t_P(x),x).   \label{eq:diffP}
\end{equation}

We need $\frac{\partial t_P}{\partial x_j}$.
 We differentiate (\ref{eq:rowTp})   to obtain
\begin{eqnarray}
  \sum_{k=1}^n \frac{\partial \alpha}{\partial x_k}(P(x))
    \left(f_k(P(x)) \frac{\partial t_P}{\partial x_j}(x) +
    \frac{\partial \varphi_k}{\partial x_j}(t_P(x),x) \right)= 0, \nonumber \\
  \left(\nabla \alpha (P(x))  \cdot f(P(x))\right)
\frac{\partial t_P}{\partial x_j}(x) +
    \sum_{k=1}^n \frac{\partial \alpha}{\partial x_k}(P(x))
                    \frac{\partial \varphi_k}{\partial x_j}(t_P(x),x) =
                    0. \label{eq:dertpuw}
\end{eqnarray}
Hence
\begin{equation}
 \frac{\partial t_P}{\partial x_j}(x) = - \frac{1}{\ip{\nabla \alpha (P(x))}{f
 (P(x))}}
 \sum_{k=1}^n  \frac{\partial \alpha}{\partial x_k}(P(x))
 \frac{\partial \varphi_k}{\partial x_j}(t_P(x),x).  \label{eq:dtpxj}
\end{equation}

\subsection{Higher order derivatives of the Poincar\'e map}
To make formulas  transparent we will drop arguments of functions
in this section, but reader should be aware that for $t_P$ and its
partial derivatives the argument is $x$, for $\varphi$ and $D_a
\varphi$ the argument is always the pair $(t_P(x),x)$.

 From \eqref{eq:diffP} we obtain
\begin{eqnarray*}
 D_{(j,c)}P &=& \frac{\partial^2}{\partial t^2}\varphi D_{(j)}t_PD_{(c)}t_P
    +\frac{\partial}{\partial t} D_{(c)}\varphi D_{(j)}t_P
    +\frac{\partial}{\partial t}\varphi D_{(j,c)}t_P\\
    &+& \frac{\partial}{\partial t}D_{(j)}\varphi D_{(c)}t_P
    +D_{(j,c)}\varphi.
\end{eqnarray*}
It is easy to see that partial derivatives of high order give rise
to  quite  complex expressions and it is not entirely obvious how
to organize it in some coherent and programmable way. For this
purpose we use the following
\begin{lemma}\label{lem:DaP}
For a multipointer $a\in\mathcal N^n_p$ we have
\begin{equation}\label{eq:DaP}
\begin{array}{rcl}
    D_aP &=& D_a\varphi + \frac{\partial \varphi}{\partial t} D_at_P\\
        &+&\sum_{k=2}^p\frac{\partial^k \varphi}{\partial t^k}\sum_{(\delta_1,\ldots,\delta_k)\in\mathcal N^p(k)}
            \prod_{j=1}^k D_{a_{\delta_j}}t_P\\
        &+&\sum_{k=2}^p\sum_{(\delta_1,\ldots,\delta_k)\in\mathcal N^p(k)}
            \sum_{s=1}^k\frac{\partial^{k-1}}{\partial
            t^{k-1}}D_{a_{\delta_s}}\varphi
            \prod_{j\neq s}D_{a_{\delta_j}}t_P
\end{array}
\end{equation}
\end{lemma}
\textbf{Proof:} By induction on $p$. For $p=1$ formula
\eqref{eq:DaP} is equivalent to (\ref{eq:diffP}), because the two
last sums are taken over empty set. Assume \eqref{eq:DaP} holds true
for some $p\geq1$ and fix $a\in\mathcal N^n_{p+1}$. Our goal is to
show that
\begin{equation*}
D_aP = R_1+R_2+R_3
\end{equation*}
where
\begin{eqnarray*}
R_1 &=& D_a\varphi+\frac{\partial}{\partial t}\varphi D_at_P\\
R_2 &=& \sum_{k=2}^{p+1} \frac{\partial^k}{\partial t^k}\varphi
    \sum_{(\delta_1,\ldots,\delta_k)\in\mathcal N^{p+1}(k)}
        \prod_{j=1}^kD_{a_{\delta_j}}t_P\\
R_3&=& \sum_{k=2}^{p+1}\sum_{(\delta_1,\ldots,\delta_k)\in\mathcal
    N^{p+1}(k)}
        \sum_{s=1}^k\frac{\partial^{k-1}}{\partial
        t^{k-1}}D_{a_{\delta_s}}\varphi\prod_{j\neq
        s}D_{a_{\delta_j}}t_P
\end{eqnarray*}

Write $a=\beta+\gamma$, where $\beta\in\mathcal N^n_p$ and
$\gamma=(a_{p+1})\in\mathcal N^n_1$. From the induction assumption
we have
\begin{equation*}
\begin{array}{rcl}
    D_aP &=& D_\gamma\left(D_\beta\varphi + \frac{\partial}{\partial t}\varphi D_\beta t_P\right)\\
        &+&D_\gamma\left(
            \sum_{k=2}^p \frac{\partial^k}{\partial t^k}\varphi
            \sum_{(\delta_1,\ldots,\delta_k)\in\mathcal N^p(k)}
        \prod_{j=1}^kD_{\beta_{\delta_j}}t_P
            \right)\\
        &+&D_\gamma\left(
           \sum_{k=2}^p\sum_{(\delta_1,\ldots,\delta_k)\in\mathcal
    N^p(k)}
        \sum_{s=1}^k\frac{\partial^{k-1}}{\partial
        t^{k-1}}D_{\beta_{\delta_s}}\varphi\prod_{j\neq
        s}D_{\beta_{\delta_j}}t_P
        \right)\\
        &=& \sum_{i=1}^{10}S_i
        \end{array}
\end{equation*}
where
\begin{equation*}
\begin{array}{rcl}
      S_1 &=& D_a\varphi+\frac{\partial}{\partial t}\varphi D_at_P\\
      S_2 &=&\frac{\partial}{\partial t}D_\beta\varphi D_\gamma t_P\\
      S_3 &=&\frac{\partial^2}{\partial t^2}\varphi D_\beta t_P
      D_\gamma t_P\\
      S_4&=& \frac{\partial}{\partial t}D_\gamma\varphi D_\beta t_P\\
      S_5  &=&
        \sum_{k=2}^p \frac{\partial^k}{\partial t^k}D_\gamma\varphi
           \sum_{(\delta_1,\ldots,\delta_k)\in\mathcal
            N^p(k)}
            \prod_{j=1}^kD_{\beta_{\delta_j}}t_P\\
      S_6  &=&
        \sum_{k=2}^p \frac{\partial^{k+1}}{\partial
        t^{k+1}}\varphi D_\gamma t_P
           \sum_{(\delta_1,\ldots,\delta_k)\in\mathcal
            N^p(k)}
            \prod_{j=1}^kD_{\beta_{\delta_j}}t_P\\
       S_7 &=&
        \sum_{k=2}^p \frac{\partial^k}{\partial t^k}\varphi
           \sum_{(\delta_1,\ldots,\delta_k)\in\mathcal
            N^p(k)}\sum_{s=1}^k D_{\beta_{\delta_s}+\gamma}t_P
            \prod_{\substack{j=1\\j\neq s}}^kD_{\beta_{\delta_j}}t_P\\
       S_8 &=&
           \sum_{k=2}^{p}\sum_{(\delta_1,\ldots,\delta_k)\in\mathcal
            N^p(k)}
                \sum_{s=1}^k\frac{\partial^{k-1}}{\partial
                t^{k-1}}D_{\beta_{\delta_s}+\gamma}\varphi\prod_{j\neq
                s}D_{\beta  _{\delta_j}}t_P\\
       S_9 &=&
           \sum_{k=2}^{p}\sum_{(\delta_1,\ldots,\delta_k)\in\mathcal
            N^p(k)}
                \sum_{s=1}^k\frac{\partial^{k}}{\partial
                t^{k}}D_{\beta_{\delta_s}}\varphi D_\gamma t_P\prod_{j\neq
                s}D_{\beta  _{\delta_j}}t_P\\
       S_{10} &=&
           \sum_{k=2}^{p}\sum_{(\delta_1,\ldots,\delta_k)\in\mathcal
            N^p(k)}\\&&
                \sum_{s=1}^k\sum_{\substack{r=1\\r\neq s}}^k\frac{\partial^{k-1}}{\partial
                t^{k-1}}D_{\beta_{\delta_s}}\varphi D_{\beta_{\delta_r}+\gamma}t_P
                \prod_{\substack{j\neq s\\ j\neq r}}D_{\beta
                _{\delta_j}}t_P
\end{array}
\end{equation*}
Obviously $R_1=S_1$. We will show that $R_2=S_3+S_6+S_7$ and
$R_3=S_2+S_4+S_5+S_8+S_9+S_{10}$.

Denote by $R_{i,k}$, $i=2,3$ a part of sum $R_i$ with fixed
$k=2,\ldots,p+1$. Similarly, let us denote by $S_{i,k}$ a part of
sum $S_i$, $i=5,\ldots,10$, for $k=2,\ldots,p$.

Using decomposition of $\mathcal N^{p+1}(2)$ as in
\eqref{eq:indexDecomposition} we obtain that $R_{2,2}=S_3 +
S_{7,2}$. Similarly, using \eqref{eq:indexDecomposition} we
observe that $R_{2,k}=S_{6,k-1}+S_{7,k}$ for $k=3,\ldots,p$.
Finally, since $\mathcal N^{p+1}(p+1)=\{((1),(2),\ldots,(p+1))\}$
and $\gamma=(a_{p+1})$ we find that $R_{2,p+1}=S_{6,p}$. This
shows that $R_2=S_3+S_6+S_7$.

It remains to show that $R_3=S_2+S_4+S_5+S_8+S_9+S_{10}$. We will
classify possible terms by the fact, where  $p+1$ appears in
$\delta_i$, $i=1,\dots,k$ and how this $\delta_i$ enters in $R_3$ as
$\delta_s$ or $\delta_j$. There are four cases
\begin{enumerate}
\item[1.] $\delta_s=(p+1)$
\item[2.] $\delta_j=(p+1)$
\item[3.] $p+1 \in \delta_s$, $|\delta_s| \geq 2$
\item[4.] $p+1 \in \delta_j$, $|\delta_j| \geq 2$
\end{enumerate}
Let us fix $k=2$. Let $(\delta_1,\delta_2)\in\mathcal N^{p+1}(2)$.
The term for case 1 is $S_4$, for case 2 is $S_2$, case 3 is
$S_{8,2}$ and case 4 is $S_{10,2}$. Hence,
$R_{3,2}=S_2+S_4+S_{8,2}+S_{10,2}$.

For $k=3,\ldots,p$ and fixed $(\delta_1,\ldots,\delta_k)\in\mathcal
N^{p+1}(k)$ we have: case 1 is given by $S_{5,k-1}$, case 2 by
$S_{9,k-1}$, case 3 by $S_{8,k}$ and case 4 by $S_{10,k}$ Hence, for
$k=3,\ldots,p$ we have
$R_{3,k}=S_{5,k-1}+S_{9,k-1}+S_{8,k}+S_{10,k}$.

Finally, for $k=p+1$ we observe, that $R_{3,p+1}=S_{5,p}+S_{9,p}$.
Indeed, in this case
$(\delta_1,\ldots,\delta_{p+1})=((1),(2),\ldots,(p+1))$. Hence,
either for $\delta_s=\gamma$ we have term $S_{5,p}$ and $\delta_s
\neq\gamma$ we have $S_{9,p}$.

We have showed that $R_3=S_2+S_4+S_5+S_8+S_9+S_{10}$ and the proof
is finished. \qed

Hence, if we know all the partial derivatives of $t_P$ up order
$p$ we can compute the partial derivatives of the Poincar\'e map
up the same order. In next subsection we show how to compute
partial derivatives of $t_P$ for affine sections.

\subsection{Partial derivatives of $t_P$ for affine sections}
Assume $\alpha:\mathbb R^n\to\mathbb R$ is an affine map given by
\begin{equation*}
\alpha(x)=\alpha_0+\sum_{i=1}^n\alpha_ix_i.
\end{equation*}
This is a quite restrictive assumption about sections, but it
leads to relatively simple formulas for $D_a t_P$ and it is
sufficient for the applications we have in mind.

\begin{lemma}\label{lem:DatP}
For a multipointer $a\in\mathcal N^n_p$ holds
\begin{equation*}
    \begin{array}{rcl}
    -D_at_P\ip{\nabla\alpha}{\frac{\partial}{\partial t}\varphi} &=&
    \ip{\nabla\alpha}{D_a\varphi}\\
    &+&
    \sum_{k=2}^p \ip{\nabla\alpha}{\frac{\partial^k}{\partial t^k}\varphi}
    \sum_{(\delta_1,\ldots,\delta_k)\in\mathcal N^p(k)}
        \prod_{j=1}^kD_{a_{\delta_j}}t_P\\
    &+&
    \sum_{k=2}^{p}\sum_{(\delta_1,\ldots,\delta_k)\in\mathcal
    N^p(k)}\\&&
        \sum_{s=1}^k\ip{\nabla\alpha}{\frac{\partial^{k-1}}{\partial
        t^{k-1}}D_{a_{\delta_s}}\varphi}\prod_{j\neq
        s}D_{a_{\delta_j}}t_P
    \end{array}
\end{equation*}
\end{lemma}
\textbf{Proof:} The proof is a direct consequence of
Lemma~\ref{lem:DaP} and \eqref{eq:rowTp}. Since $\alpha$ is
affine, by differentiating of $\alpha(P(x))=C$ we get
$\ip{\nabla\alpha}{D_aP}=0$. Using  formula \eqref{eq:DaP} for
$D_aP$ we obtain our assertion. \qed

Fix $[x] \subset \mathbb R^n$ and assume we have a rigorous bound
for $t_P([x])\in[t_1,t_2]$ (see  \cite[Section 6]{Z} for more
details on this).
 Lemmas~\ref{lem:DatP} and~\ref{lem:DaP} show that given  rigorous bounds
  for the partial derivatives $D_a\varphi([t_1,t_2],[x])$ and
$\frac{\partial^k}{\partial t^k}D_a\varphi([t_1,t_2],[x])$ up to
some order $p$ we can compute recursively rigorous bounds for the
partial derivatives of $t_P([x])$ and $P([x])$ up to the same order.
Notice, that $\frac{\partial^k}{\partial t^k}D_a\varphi$ are given
by Taylor coefficients of the solution of \eqref{eq:mainEquation}
with  initial conditions $P([x])$ for $\mathcal C^0$ part and
$D_a\varphi(t_P(x),[x])$ for equations for variations. Hence, these
coefficients can be easily computed using the automatic
differentiation algorithm.

%% file: pendulum.tex
\section{Applications.}
\label{sec:pendulum}
 One of the typical invariant sets in
hamiltonian mechanics are invariant tori.  However, the existence
of invariant torus in a given system is often difficult to prove
despite the fact that the theory is quite well developed. Probably
the best work in this direction was done by Celletti and Chercia
\cite{CC1,CC2}, where the an effective application (computer
assisted proof) of KAM theory to the restricted three body problem
modelling system consisting of Sun,  Jupiter  and asteroid 12
Victoria was given. Our aim here is more modest as we focus on the
invariant tori emanating from the elliptic fixed point satisfying
suitable twist condition.

In this section we show that the rigorous computations of
 partial derivatives of a dynamical system up to order $3$ or $5$ can be used to prove
that in a particular system an invariant torus exists around some
elliptic periodic orbits. In this section this will be done for
the forced pendulum equation and the Michelson system.

\subsection{Area preserving maps on the plane, normal forms and  KAM theorem }

\begin{definition}
Let $f:\mathbb{R}^2 \to \mathbb{R}^2$ be a smooth area preserving
map, such that $f(p)=p$. Let $\lambda$ and $\mu$ be eigenvalues of
$df(p)$. Following \cite{SM} we will call the point $p$
\begin{itemize}
\item \emph{hyperbolic} if $\lambda,\mu \in \mathbb{R}$ and $\lambda\neq
\mu$,
\item \emph{elliptic} if $\lambda=\overline{\mu}$ and $\lambda\neq
\mu$,
\item \emph{parabolic} if $\lambda=\mu$.
\end{itemize}
\end{definition}

The following KAM theorem will be the main tool to prove the
existence of invariant tori in this paper.

\begin{theorem}\cite[\S 32]{SM}\label{thm:SM1}
Consider an analytic area preserving map $f:\mathbb R^2\to\mathbb
R^2$,  $f(r,s)=(r_1,s_1)$ where
\begin{eqnarray}
r_1 &=& r\cos \alpha - s\sin \alpha + O_{2l+2}  \nonumber \\
s_1 &=& r\sin \alpha + s\cos \alpha + O_{2l+2} \label{eq:normalForm}  \\
\alpha &=& \sum_{k=0}^l\gamma_k\left(r^2+s^2\right)^k \nonumber
\end{eqnarray}
and $O_{2l+2}$ denotes convergent power series in $r,s$ with terms
of order greater than $2l+1$, only.

If at least one of $\gamma_1,\ldots,\gamma_l$ is not zero then the
origin is a stable fixed point for map $f$. Moreover, in any
neighborhood $U$ of point $0$ there exists an invariant curve for
map $f$ around the origin contained in $U$.
\end{theorem}

The next theorem and its proof tells how to bring a planar area
preserving map in the neighborhood of an elliptic fixed point into
the form (\ref{eq:normalForm}).
\begin{theorem}\cite[\S 23]{SM}\label{thm:SM2}
Consider an analytic area preserving map $f:\mathbb R^2\to\mathbb
R^2$ such that $f(0)=0$. Let $\lambda,\bar\lambda$ be complex
eigenvalues of $Df(0)$, such that $|\lambda|=|\bar\lambda|=1$. If
$\lambda^k\neq 1$ for $k=1,\ldots,2l+2$, then there is an analytic
area preserving substitution such that in the new coordinates
mapping $f$ has form \eqref{eq:normalForm}.
\end{theorem}
The proof of the above theorem is constructive, i.e. given the
power series for $f$ at an elliptic fixed point one can construct
explicitly an area preserving substitution and compute the
coefficients $\gamma_0,\ldots,\gamma_l$ in \eqref{eq:normalForm}.
An explicit formula for the coefficient $\gamma_1$ in the above
normal form is given in  Appendix \ref{sec:appendixA}.

\subsection{The existence of invariant tori in forced pendulum.}
Consider an equation
\begin{equation}\label{eq:pendulum}
\ddot\theta =  -\sin(\theta)+\sin(\omega t)
\end{equation}
Observe that (\ref{eq:pendulum}) is hamiltonian.

Let us denote by $P_\omega:\mathbb R^2\oarrow\mathbb R^2$ the
Poincar\'e map for Equation \eqref{eq:pendulum} with a parameter
$\omega$, i.e. $P_\omega=\varphi(2\pi/\omega,\cdot)$, where
$\varphi:\mathbb R\times\mathbb R^2\oarrow\mathbb R^2$ is a local
flow induced by \eqref{eq:pendulum}. Observe that
(\ref{eq:pendulum}) is nonautonomous, but it is equivalent to
first order system of autonomous ODE given by
\begin{eqnarray}
   \frac{d\theta}{ds}&=&v  \nonumber \\
   \frac{dv}{ds}&=& - \sin(\theta) + \sin(\omega t)  \label{eq:pendsys} \\
   \frac{dt}{ds}&=& 1. \nonumber
\end{eqnarray}
In the sequel all rigorous computations for (\ref{eq:pendulum})
will be in fact performed for the system (\ref{eq:pendsys}).

Observe that to any invariant closed curve for $P_\omega$
corresponds and invariant 2-torus for (\ref{eq:pendulum}).

Consider a set of parameter values
\begin{eqnarray*}
 \Omega_1 &=&[2,2.994],\quad
 \Omega_2=[3,3.997],\quad\Omega_3=[4,8] \\
  \Omega&=&\Omega_1\cup\Omega_2\cup\Omega_3
\end{eqnarray*}

The following lemma was proved with computer assistance
\begin{lemma}\label{lem:pendulum}
For all parameter values $\omega\in\Omega$ there exists an elliptic
fixed point $x_\omega\in\mathbb R^2$ for $P_\omega$. Moreover, there
exists an area-preserving substitution such that in the new
coordinates the map $f_\omega(x)=P_\omega(x+x_\omega)-x_\omega$ has
the form \eqref{eq:normalForm} with $l=1$ and $\gamma_1\neq0$.
\end{lemma}
Before we give the proof, let us briefly comment about the choice of
the parameter set $\Omega$. For parameter values slightly lower than
$2$ we observe the parabolic case, i.e. there exists a parameter
value $\omega_1$ for which eigenvalues of the derivative of
$P_{\omega_1}$ are equal to $-1$. In two gaps in $\Omega$ below $3$
and $4$ we have resonances of low order. Namely,  we have parameter
values with an elliptic fixed with eigenvalues  to
$e^{\pm2\pi/3}=\frac{-1}{2}\pm\frac{\sqrt{3}}{2}i$ and $e^{\pm
i\pi/2}=\pm i$, respectively. Clearly, in a computer assisted proof
we need to exclude a small interval around those parameters. For
$\omega
>4$ it seems that the interval $\Omega_3$ can be extended much
further to the right without any difficulty.

\noindent \textbf{Proof of Lemma~\ref{lem:pendulum}:} A computer
assisted proof consists of the following steps. We cover the set
$\Omega$ by $9910$ nonequal subintervals $\omega_i$. Diameters of
$\omega_i$'s were relatively large for values  far away from the
parabolic cases and very small close to them. For a fixed
subinterval $\omega_i$ we proceed as follows
\begin{enumerate}
\item Let $\bar\omega$ denote an approximate center of the interval $\omega_i$.
We find an approximate fixed point for $P_{\bar\omega}$ using the
standard nonrigorous Newton method. Let us denote such a point by
$x_i$.
\item We define a box centered at $x_i$, i.e we set
$v_i:=x_i+[-\varepsilon_i,\varepsilon_i]^2$, where
$\varepsilon_i>0$ depends on subinterval $\omega_i$ - the values
we used are from the interval $[5\cdot10^{-6},3\cdot10^{-3}]$,
depending on whether $x_i$ close to parameter values corresponding
to parabolic cases.
\item Using the $\mathcal C^1$-Lohner algorithm we compute
the Interval Newton operator \cite{Mo,N,A}
$N_i:=N(P_{\omega_i}-\Id,x_i,v_i)$ and verify that $N_i\subset
\inte v_i$. This proves that for all $\omega\in\omega_i$ there
exists a unique fixed point $x_\omega\in N_i$ for $P_\omega$.
\item Using the $\mathcal C^3$-Lohner algorithm we compute
a rigorous bound for $P_{\omega_i}(N_i)$ and $D^\alpha
P_{\omega_i}(N_i)$, $\alpha\in\mathbb N^2_1\cup\mathbb
N^2_2\cup\mathbb N^2_3$. Hence, we obtain a rigorous bound for the
coefficients in
\begin{equation*}
f_\omega(x) = \sum_{\substack{|\alpha|=1\\
\alpha\in\mathbb N^2}}^3 \frac{1}{\alpha!}D_a P(x_\omega)x^ + O_{4}
\end{equation*}
\item We show that an arbitrary matrix $M\in DP_{\omega_i}(N_i)$ has
a pair of complex eigenvalues $\lambda, \bar\lambda$ which satisfy
$\lambda^k\neq1$ for $k=1,\ldots,4$. From Theorem \ref{thm:SM2} it
follows there exists an area-preserving substitution such that in
the new coordinates the map $f_\omega$ for $\omega\in\omega_i$ has
the form \eqref{eq:normalForm} with $l=1$.
\item We compute a rigorous bound for $\gamma_0$ and $\gamma_1$
which appear in the formula \eqref{eq:normalForm} and verify that
for $\omega\in\omega_i$ holds $\gamma_1\neq0$.
\end{enumerate}
The rigorous bounds for the values of $\gamma_1$ on $\Omega$ are
\begin{eqnarray*}
\gamma_1(\Omega_1)&\subset&
[0.29930416771330087,30.118260918229566]\\
\gamma_1(\Omega_2)&\subset&
[0.099747909112924596,0.56550301088840627]\\
\gamma_1(\Omega_3)&\subset&[0.18574835001593507,0.4129279974577012]
\end{eqnarray*}
 A computer assisted proof of the above took approximately $95$
minutes on the Pentium IV 3GHz processor. \qed

As a straightforward consequence of Lemma~\ref{lem:pendulum} and
Theorem \ref{thm:SM1} we obtain
\begin{theorem}
For all parameter values $\omega\in\Omega$ there exists an
elliptic fixed point $x_\omega\in\mathbb R^2$ for $P_\omega$.
Moreover, any neighborhood of point $x_\omega$ contains
 an invariant curve for $P_\omega$ around $x_\omega$.
\end{theorem}

\subsection{Higher order normal forms.}
In the previous section it was shown  that $\mathcal C^3$
computations are  sufficient to prove that for (\ref{eq:pendulum})
a family of invariant tori exists. However, it may happen that the
coefficient $\gamma_1$ in the normal form vanishes. In this
situation we may try to compute higher order normal form. As an
example we consider a pendulum with a different forcing term,
\begin{equation}\label{eq:pendulum2}
\ddot \theta = -\sin(\theta) + \sin(\omega t) + \sin(2\omega t).
\end{equation}
\begin{theorem}
Let $P_\omega$ be the Poincar\'e map for \eqref{eq:pendulum2}. For
all parameter values
$\omega\in\Omega_*=[2.9957694795,2.9957694796]$ there exists an
elliptic fixed point $x_\omega\in\mathbb R^2$ for $P_\omega$.
Moreover, any neighbourhood of point $x_\omega$ contains an
invariant curve for $P_\omega$ around $x_\omega$.
\end{theorem}
\textbf{Proof:} The main concept of the proof  is the same as in
Lemma~\ref{lem:pendulum}. Using the nonrigorous Newton method we
find an approximate fixed point
\begin{equation*}
x=(-7.7491573604896152\cdot10^{-12},-0.54723831527031352).
\end{equation*}
We set $v=x+3\cdot10^{-5}([-1,1]\times[-1,1])$. Using the
$\mathcal C^1$-Lohner algorithm we compute the Interval Newton
Operator of $P_\omega-\Id$ on $v$ and we obtain that for all
$\omega\in\Omega_*$,
$N=N(P_{\omega}-\Id,\mathrm{center}(v),v)\subset(N_1,N_2)$, where
\begin{eqnarray*}
    N_1&=&[-5.1582932672798325,5.1582631625020222]\cdot 10^{-10}\\
    N_2&=&[-0.54723831580217108,-0.54723831470891193]
\end{eqnarray*}
Since $N\subset v$ we conclude that for all $\omega\in\Omega_*$
there exists a unique fixed point $x_\omega\in N$ for the Poincar\'e
map.

Using $\mathcal C^5$-Lohner algorithm we compute a rigorous bound
for $P_{\Omega_*}(N)$ and $D^\alpha P_{\Omega_*}(N)$,
$\alpha\in\mathbb N^2_1\cup\ldots\cup\mathbb N^2_5$. Hence, we
obtain a rigorous bound for the coefficients in
\begin{equation*}
f_\omega(x) = \sum_{\substack{|\alpha|=1\\
\alpha\in\mathbb N^2}}^5 \frac{1}{\alpha!}D^\alpha
P(x_\omega)x^\alpha + O_{6}
\end{equation*}

We show that an arbitrary matrix $M\in DP_{\Omega_*}(N)$ has a pair
of complex eigenvalues $\lambda, \bar\lambda$ which satisfy
$\lambda^k\neq1$ for $k=1,\ldots,6$. From Theorem \ref{thm:SM2} it
follows there exists an area-preserving substitution such that in
the new coordinates the map $f_\omega$ for $\omega\in\Omega_*$ has
the form \eqref{eq:normalForm} with $l=2$.

Next, we compute a rigorous bound for $\gamma_1$ and
 $\gamma_2$ which appear in the formula \eqref{eq:normalForm} and we
 get
 \begin{eqnarray*}
    \gamma_1(\Omega_*)&\subset&[-5.3924276719042241,5.381714805052106]\cdot10^{-6}\\
    \gamma_2(\Omega_*)&\subset&[199.95180660157078,199.99104965939162]
 \end{eqnarray*}
Since for $\omega\in\Omega_*$, $\gamma_2(\omega)\neq0$ the assertion
follows from Theorem~\ref{thm:SM1}.
 \qed

The main observation which makes this example interesting is that
there exists $\omega_*\in\Omega_*$ for which
$\gamma_1(\omega_*)=0$ and we cannot conclude the existence of
invariant tori for all $\omega \in \Omega_*$ from $\mathcal C^3$
computations. To be more precise, we computed the coefficient
$\gamma_1$ for the parameter values $\omega_1=\min\Omega_*$ and
$\omega_2=\max\Omega_*$ and we get
\begin{eqnarray*}
    \gamma_1(\omega_1)&\in&[-2.3559594437885885,-1.3593457220363871]\cdot10^{-8}\\
    \gamma_1(\omega_2)&\in&[2.9671154858524365\cdot10^{-9},1.2819312939263052\cdot10^{-8}]
\end{eqnarray*}
Since $\gamma_1$ exists for all $\omega\in\Omega_*$ and depends
continuously on $\omega$ we conclude, that $\gamma_1(\omega_*)=0$
for some $\omega_*\in\Omega_*$.

%% file: msystem.tex
\subsection{Application to the Michelson system}

The existence of an invariant curve for a planar map $f:\mathbb
R^2\to\mathbb R^2$ can be proven without assumption that $f$ is
measure preserving. The key assumption in the proof given in
\cite{SM} is that any curve $\gamma$ around an elliptic point
intersect its image under $f$, i.e.
$f(\gamma)\cap\gamma\neq\emptyset$. Such a situation is also
observed in reversible planar map around an symmetric
elliptic fixed points.

\begin{definition}
An invertible transformation $M:\Omega\longrightarrow \Omega$ is
called a \emph{reversing symmetry} of a local dynamical system
$\phi:\mathbb T\times\Omega\longrightarrow\Omega$, $\mathbb
T=\mathbb R$ or $\mathbb T= \mathbb Z$ if the following conditions
are satisfied
\begin{enumerate}
\item if $(t,x)\in \dom(\phi)$ then $(-t,S(x))\in\dom(\phi)$.
\item $S(\phi(t,x)) = \phi(-t,S(x)))$
\end{enumerate}
\end{definition}
\begin{remark}
In the discrete time case, the above two conditions are equivalent
to identity
\[
M\circ f=f^{-1}\circ M.
\]
where $f=\phi(1,\cdot)$ is a generator of $\phi$.
\end{remark}
\begin{definition}
Let $\phi:\mathbb T\times\Omega\to\Omega$ be a local (discrete or
continuous) dynamical system. For $x\in\Omega$ put
\begin{eqnarray*}
I(x)&=&\{t\in\mathbb T : (t,x)\in\dom(\phi)\}\\
\orbit{x}&=&\{\phi(t,x)\in\Omega : t\in I(x)\}
\end{eqnarray*}
The set $\orbit{x}$ will be called a \emph{trajectory} of a point
$x$.
\end{definition}
\begin{definition}
Assume $S$ is an reversing symmetry for $\phi:\mathbb
T\times\Omega\to\Omega$. An orbit $\orbit{x}$ is called
$S$-symmetric orbit if $\orbit{x}=S(x)$.
\end{definition}
\begin{remark}\cite{La}
In continuous case the orbit $\orbit{x}$ is $S$-symmetric if it
contains a point from the set $\Fix(S)=\{y:S(y)=y\}$.
\end{remark}

\begin{remark}\cite[Lem.3.3]{Wi2}\label{rem:revPoincare}
It is easy to see that if $\Theta\subset\Omega$ is a Poincar\'e
section for a $R$-reversible flow $\phi:\mathbb
R\times\Omega\to\Omega$ such that $\Theta=R(\Theta)$ then the
Poincar\'e map $P:\Theta\to\Theta$ is $R|_{\Theta}$-reversible.
\end{remark}

As we observed at the beginning of this section, an $R$-reversible
planar map may admit an invariant curve around an $R$-symmetric
elliptic fixed point. In reversible case a planar map admits the
same normal form around symmetric, elliptic fixed point as in the
area-preserving case and the substitution which tends the map to
the normal form is exactly the same as we described in
Appendix~\ref{sec:appendixA} -- for details see \cite{Se,BHS}.

Consider an ODE
\begin{equation}\label{eq:michelson}
\left\{
    \begin{array}{rcl}
        \dot x &=& y\\
        \dot y &=& z\\
        \dot z &=& c^2-y-\frac{1}{2}x^2
    \end{array}
\right.
\end{equation}
On one hand, the system \eqref{eq:michelson} is an
equation for the steady state solution of one-dimensional
Kuramoto-Sivashinsky PDE and it is known in the literature as the
Michelson system\cite{Mi}. On the other hand, this system
appears as a part of the limit family of the unfolding of the
nilpotent singularity of codimension three (see \cite{DIK1}).

The system \eqref{eq:michelson} is reversible with respect to the
symmetry
\begin{equation}\label{eq:michelsonSymmetry}
    R:(x,y,z,t)\to(-x,y,-z,-t)
\end{equation}
and since the divergence vanishes it is also volume preserving.

A dynamical system induced by \eqref{eq:michelson} exhibits several
types of dynamics for different values of parameter. For
sufficiently large $c$ there is a simple invariant set consisting of
two equilibria $(\pm c\sqrt2,0,0)$ and heteroclinic orbit between
them \cite{MC}. Lau \cite{Lau} numerically observed that when the
parameter $c$ decreases a cascade of cocoon bifurcations occurs and
at the limit value $c\approx 1.266232337$ a periodic orbit is born
through a saddle-node bifurcation. This hypothesis has been proved
in \cite{KWZ}. The computer assisted proof of this fact given in
\cite{KWZ} uses the algorithm presented in this paper in order to
compute partial derivatives up to second order for a certain
Poincar\'e map.

For the parameter value equal to one and slightly smaller than one
it was proven in \cite{DIK,Wi1,Wi2,Wi3} that the system has rich
and complicated dynamics including symbolic dynamics, heteroclinic
solutions, Shilnikov homoclinic solutions.

However, as the bifurcations diagram presented by Michelson
suggests \cite[Fig.1]{Mi} for all parameter values
$c\in(0,0.3195)$ there are at least two elliptic periodic orbits
with large invariant islands around them. In this section we
present a proof that such islands exist for some range of
parameter values. The main idea of the proof is almost the same as
in the previous section. There are two main differences. First,
the Poincar\'e map will not be a time shift. Therefore
computations of the partial derivatives of the Poincar\'e map
require Lemma~\ref{lem:DaP} and Lemma~\ref{lem:DatP}. Second
difference is: we use the shooting method instead of the interval
Newton method for the proof of the existence of symmetric periodic
orbit.

The aim of this section is to prove the following
\begin{theorem}\label{thm:toriMS}
For all parameter values from the set
\begin{equation*}
C=C_1\cup C_2 = [0.1,0.225]\cup[0.226,0.25]
\end{equation*}
there exists a symmetric elliptic periodic orbit for the Michelson
system \eqref{eq:michelson}. Moreover, each neighbourhood of such an
orbit contains a $2D$ tori invariant under the flow generated by the
Michelson system.
\end{theorem}

Let us define the Poincar\'e section $\Pi:=\{(0,z,y):z,y\in\mathbb
R\}$. Let $P_c=(P_1,P_2):\Pi\oarrow\Pi$ be the Poincar\'e map for
the system with the parameter value $c$. Notice, that $P_c$ is in
fact a half Poincar\'e map, which means that the trajectory of $x$
crosses $\Pi$ in opposite directions when passing through $x$ and
$P_c(x)$, and therefore periodic orbits  for the Michelson system
corresponds to periodic points for $P^2$.

Since the section $\Pi$ is invariant under symmetry
$(x,y,z)\to(-x,y,-z)$, from Remark~\ref{rem:revPoincare} the
Poincar\'e map is also reversible with respect to an involution
$R(y,z)=(y,-z)$. We will use the same letter $R$ to denote the
reversing symmetry of the Poincar\'e map and the Michelson system.

\begin{figure}
    \centerline{\includegraphics{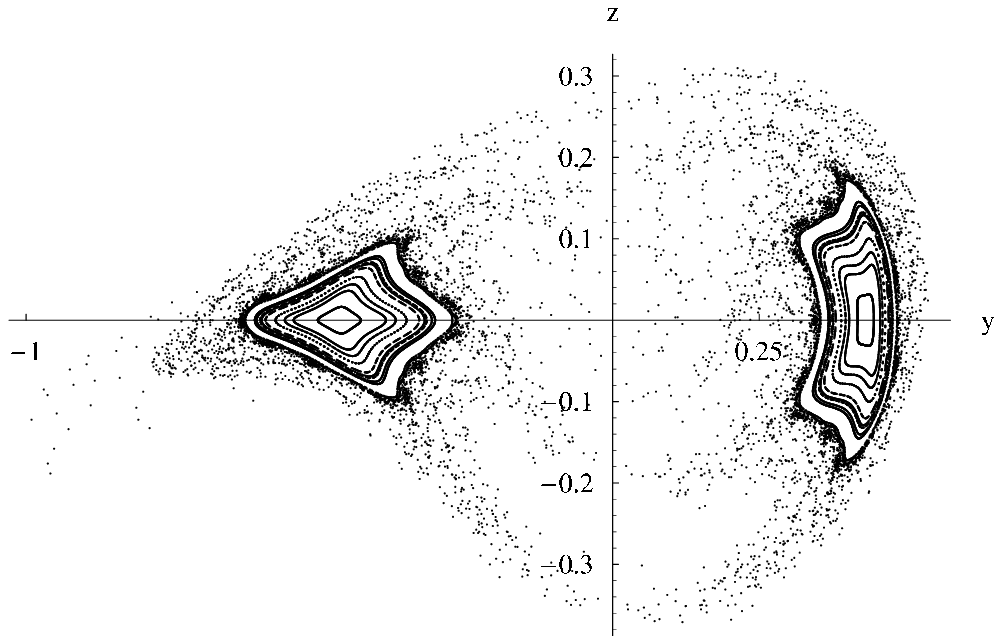}}
    \centerline{\includegraphics{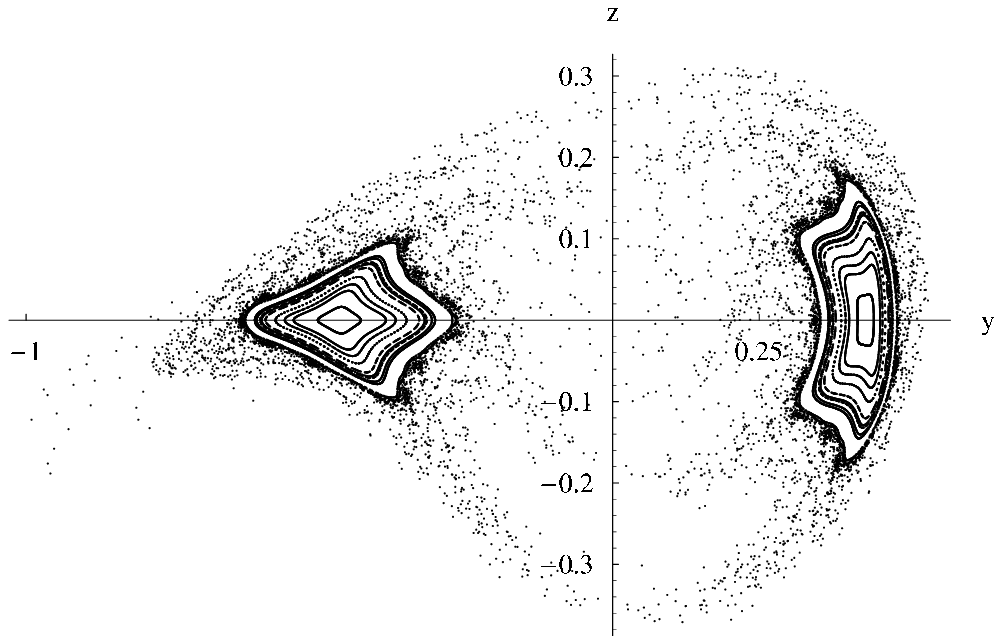}}
    \caption{
        Phase portrait of the Poincar\'e map $P_c$ (top) before bifurcation for $c=0.225$ and (bottom)
        after bifurcation for $c=0.226$ with four periodic islands.
        Between those parameters resonant case
        occurs with eigenvalues equal to $\pm i$.
        See also auxiliary material \cite{Wi4}.
        \label{fig:mpp}
    }
\end{figure}

Let us comment about the choice of the set $C$. In the gap between
intervals $C_1$ and $C_2$ there is a parameter value $c_*$ for which
the eigenvalues of the Poincar\'e map $P^2_{c_*}$ are $\pm i$.
Apparently at this parameter value we have a bifurcation and four
periodic islands are born as it is shown in Fig.\ref{fig:mpp} - see
also a movie {\tt mpp.mov} available at \cite{Wi4} which presents an
animation of the phase portrait of $P_c$ for the parameter values
from the range $[0.1,0.25]$.

\noindent\textbf{Proof of Theorem~\ref{thm:toriMS}.} The main
concept of the proof is quite similar to the one presented in
Lemma~\ref{lem:pendulum}. We divide the set $C$ of parameter values
onto $20800$ nonequal parts (smaller when close to the bifurcation
parameter $c_*$ and close to $0.1$ and $0.25$). For a fixed
subinterval $c_i$ from the grid we proceed as follows
\begin{enumerate}
\item Let $\bar c$ denote a center of the interval $c_i$. We find an approximate fixed point
of $P^2_{\bar c}$ using the standard nonrigorous Newton method. Let
us denote this point by $(y_i,z_i)$.
\item Since the map $P_c$ is reversible one can prove the existence
of the fixed point for $P_c^2$ using the shooting method as follows.

 Let
$\mathrm{Fix}(R)=\{(y,z)\in\Pi:R(y,z)=(y,z)\} =
\{(y,0)\in\Pi:y\in\mathbb R\}$. Since $P_c$ satisfies $(P_c\circ
R)^2=\Id$ whenever the left side is defined, one can see that if
$x\in\mathrm{Fix}(R)$ and $P_c(x)\in\mathrm{Fix}(R)$ then
$P^2_c(x)=x$. Let us remark, that we always get an approximate fixed
points $(y_i,z_i)$ resulting from the nonrigorous Newton method very
close to $\mathrm{Fix}(R)$. We define two points
$u_1=(y_i-\varepsilon_i,0),
u_2=(y_i+\varepsilon_i,0)\in\mathrm{Fix}(R)$, where $\varepsilon_i$
is a small number depending on $c_i$ and we show that
$\pi_z(P_{c_i}(u_1))\cdot\pi_z(P_{c_i}(u_2))<0$, where $\pi_z$ is a
projection onto $z$ coordinate. Hence, if the $P_{c_i}$ is defined
on the set $N_i=(0,[y_i-\varepsilon_i,y_i+\varepsilon_i],0)$ then
for all parameter values $c\in c_i$ there is a point $u_c\in N$
which satisfies $\pi_z(P_c(u_c))=0$ and therefore
$P_c(u_c)\in\mathrm{Fix}(R)$. This shows that for all $c\in c_i$
there exists a fixed point for $P^2_c$ inside $N_i$ provided $P_c$
is defined on $N_i$, which will be discussed below.
\item Using $\mathcal C^3$-Lohner algorithm we compute rigorous bounds for
$P^2_{c_i}(N_i)$ and $D^\alpha P^2_{c_i}(N_i)$ for $\alpha\in\mathbb
N^2_1\cup\mathbb N^2_2\cup\mathbb N^2_3$. This implies also that
$N_i \subset \dom P_{c_i}$.
\item We show that an arbitrary matrix $M\in DP_{c_i}(N_i)$ has
a pair of complex eigenvalues $\lambda, \bar\lambda$ which satisfy
$\lambda^k\neq1$ for $k=1,\ldots,4$. From Theorem \ref{thm:SM2} it
follows there exists an area-preserving substitution such that in
the new coordinates the map $P_c$ for $c\in c_i$ has the form
\eqref{eq:normalForm} with $l=1$.
\item We compute a rigorous bound for $\gamma_0$ and $\gamma_1$
which appear in the formula \eqref{eq:normalForm} and verify that
for $c\in c_i$ holds $\gamma_1\neq0$.
\end{enumerate}
The rigorous bounds for the values of $\gamma_1$ on $C$ are
\begin{eqnarray*}
\gamma_1(C_1)&\subset& [0.014515898754816965,157.76639522562903]\\
\gamma_1(C_2)&\subset& [1.1002393483255526,151.35147664498677]
\end{eqnarray*}
The computer assisted proof of the above took approximately $7$
hours and $50$ minutes on the Pentium IV 3GHz processor. \qed

%% file: implementation.tex
\section{Implementation notes.}
All the algorithms presented in this paper have been implemented
in {\tt C++} by authors and are part of the CAPD library
\cite{CAPD}. In particular, the package implements the computation
of partial derivatives of a flow with respect to initial
condition, partial derivatives of Poincar\'e maps for linear
sections and computations of normal forms for planar maps up to
order $5$.

The implementation combines the automatic and symbolic
differentiation in order to generate a coefficients in Taylor
series for the solutions of the system \eqref{eq:mainEquation}.

Our tests shows that without difficulty we can compute partial
derivatives up to order $3$ for an equation in $8$-dimensional
phase space (which gives $1320$ equations to solve) on a computer
with 512MB memory. However, our current implementation is
optimized for lower dimensional problems. All the trees which
represent formulas \eqref{eq:mainEquation} are stored in the
memory of a computer. This speeds up computations because we do
not need to recompute all the multiindices, multipointers and
submultipointers in each step of the algorithm. Unfortunately,
such an implementation is memory-consuming. Therefore, higher
dimensional problems require a computer with huge memory even for
$\mathcal C^3$ or $\mathcal C^5$ computations.

%% file: appendix.tex
\appendix
\section{Explicit formulas for third order normal forms for a planar
map}\label{sec:appendixA}

The goal of this section is to give some details about the proof of
Theorem~\ref{thm:SM2}. We want to present some formulas to give the
reader the feeling about the necessary computations.

Throughout this section we assume that the assumptions of
Theorem~\ref{thm:SM2} are satisfied. In the neighbourhood of $0$ $f$
is given by a real, convergent power series
\begin{eqnarray*}
f(x,y)&=&(x_1,y_1)\\
x_1&=&\sum_{k=1}^\infty\sum_{l=0}^ka_{k-l,l}x^ly^{k-l}\\
y_1&=&\sum_{k=1}^\infty\sum_{l=0}^kb_{k-l,l}x^ly^{k-l}
\end{eqnarray*}

Denote also by $f:\mathbb C^2\to\mathbb C^2$ a complex extension
of $f$. Let $\lambda,\bar\lambda\in\mathbb C$ be complex
eigenvalues of $Df(0)$ and $v,\bar v\in\mathbb C$ corresponding
eigenvectors (here bar denotes the complex conjugation). Then,
using a linear substitution of the form $L=[v^T,\bar v^T]$, we can
change the coordinate system such that in the new coordinates the
mapping $f$ has the form
\begin{eqnarray*}
    f(\xi,\eta)&=&(\lambda\xi + p(\xi,\eta),\bar\lambda\eta+q(\xi,\eta))\\
    p(\xi,\eta)&=&\sum_{k=2}^\infty\sum_{l=0}^k p_{l,k-l}\xi^l\eta^{k-l}\\
    q(\xi,\eta)&=&\sum_{k=2}^\infty\sum_{l=0}^k
    q_{l,k-l}\xi^l\eta^{k-l} \\
    \overline{p_{i,j}} &=& q_{j,i} \quad \mbox{ for $i,j\geq0$}.
    \label{eq:realitycond}
\end{eqnarray*}
The last condition is a consequence of the invariance of
$\mathbb{R}^2 \subset \mathbb{C}^2$ under the complex map $f$. We
will refer to it as  \emph{the reality condition}. Namely,  the
set $\mathbb{R}^2 \subset \mathbb{C}^2$ in the new coordinates
$(\xi,\eta)$ is given by $\xi=\overline{\eta}$ and the condition
$f(\mathbb{R}^2)\subset \mathbb{R}^2$ expressed in coordinates
$(\xi,\eta)$ is equivalent to (\ref{eq:realitycond}).

Assume now, that $\lambda^k\neq1$ for $k=1,\ldots,4$. Then an
analytic area-preserving substitution satisfying reality
condition (\ref{eq:realitycond})
\begin{eqnarray*}
    (\Phi(z,v),\Psi(z,v))&=&(z_1,v_1)\\
    z_1&=&z+\sum_{k=1}^3\sum_{l=0}^k\phi_{l,k-l}z^lv^{k-l}+\cdots\\
    v_1&=&v+\sum_{k=1}^3\sum_{l=0}^k\psi_{l,k-l}z^lv^{k-l}+\cdots
\end{eqnarray*}
where
\begin{eqnarray*}
    \overline{\psi_{2,0}}=\phi_{0,2}&=&-\lambda^2p_{0,2}(\lambda^3-1)^{-1}\\
    \overline{\psi_{1,1}}=\phi_{1,1}&=&-p_{1,1}(\lambda-1)^{-1}\\
    \overline{\psi_{0,2}}=\phi_{2,0}&=&p_{2,0}(\lambda^2-\lambda)^{-1}\\
    \overline{\psi_{3,0}}=\phi_{0,3} &=& -\lambda^3\left(
        p_{0,3}+p_{1,1}\phi_{0,2}+2q_{0,2}\psi_{0,2}
    \right)(\lambda^4-1)^{-1}\\
    \overline{\psi_{2,1}}=\phi_{1,2} &=& \frac{-\lambda}{\lambda^2-1}
    \left(
        p_{1,2}+2p_{2,0}\phi_{0,2}+p_{1,1}\phi_{1, 1}+p_{1,1}\psi_{0,2}+2p_{0,2}\psi_{1,1}
    \right)\\
    \overline{\psi_{1,2}}=\phi_{2,1} &=& -\phi_{2,0}\psi_{0,2} + \phi_{0, 2}\psi_{2, 0}\\
    \overline{\psi_{0,3}}=\phi_{3,0} &=& \left(
        p_{3,0} + 2 p_{2,0}\phi_{2,0} +p_{1,1}\psi_{2,0}
    \right)(\lambda^3-\lambda)^{-1}
\end{eqnarray*}
brings $f=(f_1,f_2)$ to the normal form
\begin{equation*}
(z,v)\to\left(z(\alpha_0+\alpha_2zv),v(\beta_0+\beta_2zv)\right)+O((zv)^2)
\end{equation*}
with
\begin{equation*}
    \begin{aligned}
    \overline{\beta_0}=\alpha_0 &= \lambda\\
    \overline{\beta_2}=\alpha_2 &= q_{1,2}+2q_{2,0}\phi_{0,2}
        +q_{1,1}\phi_{1,1}+q_{1,1}\psi_{0,2}
        +2q_{0,2}\psi_{1,1}
    \end{aligned}
\end{equation*}
Finally, let $\gamma_0\in\mathbb R$ be such that
$\lambda=\alpha_0=e^{i\gamma_0}$ and we compute coefficient
$\gamma_1$ by
\begin{equation*}
\gamma_1 = \frac{- i\alpha_2}{\alpha_0} = \frac{i\beta_2}{\beta_0}
\end{equation*}
From the proof given in \cite{SM} it follows that
$\gamma_1\in\mathbb R$ and the mapping $f$ in coordinates $(z,v)$
has the form
\begin{equation*}
    f(z,v)=\left(
        ze^{i(\gamma_0+\gamma_1zv)},
        v e^{-i(\gamma_0+\gamma_1zv)}
    \right) + O_4
\end{equation*}
where $O_4$ is a convergent power series with the terms of degree
at least $4$.

Again, the coefficients of $f(z,v)$ satisfy reality condition (\ref{eq:realitycond}). In order to express
this normal form in terms of real variables we make a linear
substitution
\begin{equation*}
    z=r+is,\qquad v = r-is
\end{equation*}
and we obtain the normal form for $f$
\begin{eqnarray*}
    f(r,s) &=& (r_1,s_1)+O_4\\
    r_1&=&r\cos(\gamma_0+\gamma_1(r^2+s^2))-s\sin(\gamma_0+\gamma_1(r^2+s^2))\\
    s_1&=& r\sin(\gamma_0+\gamma_1(r^2+s^2))+s\cos(\gamma_0+\gamma_1(r^2+s^2))
\end{eqnarray*}
which agrees with \eqref{eq:normalForm}.

The  formulas  for higher order terms  $\phi_{i,j},\psi_{i,j}$
(and for $\gamma_2$, which are not given here) has been computed
in Mathematica.